\documentclass{article}

\usepackage{amsmath, amsthm, amssymb,graphics,graphicx}

\usepackage{authblk, float, bm, caption, subcaption}

\usepackage{hyperref}
\usepackage{relsize}

\title{Poly-Sinc Solution of Stochastic Elliptic Differential Equations}

\author{Maha Youssef\thanks{maha.youssef@uni-greifswald.de: Corresponding author}{ }}
\author{Roland Pulch\thanks{roland.pulch@uni-greifswald.de}}
\affil{\small{Institute of Mathematics and Computer Science, University of Greifswald, Walther-Rathenau-Stra{\ss}e 47, 17489 Greifswald, Germany}}

\date{}

\newtheorem{thm}{Theorem}
\newtheorem*{prf}{Proof}

\newtheorem{expt}{Experiment}

\begin{document}

\maketitle
\begin{abstract}
In this paper, we introduce a numerical solution of a stochastic partial differential equation (SPDE) of elliptic type using polynomial chaos along side with polynomial approximation at Sinc points. These Sinc points are defined by a conformal map and when mixed with the polynomial interpolation, it yields an accurate approximation. The first step to solve SPDE is to use stochastic Galerkin method in conjunction with polynomial chaos, which implies a system of deterministic partial differential equations to be solved. The main difficulty is the higher dimensionality of the resulting system of partial differential equations. The idea here is to solve this system using a small number of collocation points. Two examples are presented, mainly using Legendre polynomials for stochastic variables. These examples illustrate that we require to sample at few points to get a representation of a model that is sufficiently accurate.


\end{abstract}

\maketitle


\section{Introduction}
\label{sec:1}

In many applications the values of the parameters of the problem are not exactly known. These uncertainties
inherent in the model yield uncertainties in the results of numerical simulations. Stochastic methods are one way to model these uncertainties and shall model this by random fields \cite{Adler_1981}. If the physical system is described by a partial differential equation (PDE), then the combination with the stochastic model results in a
stochastic partial differential equation (SPDE). The solution of the SPDE is again a random field, describing both the expected response and quantifying its uncertainty. SPDEs can be interpreted mathematically in several ways.

In the numerical framework, the stochastic regularity of the solution determines the convergence rate of numerical approximations, and a variational theory for this was earlier devised in \cite{Theting_2000}. The ultimate goal in the solution of SPDEs is usually the computation of response statistics, i.e. a functional of the solution. Monte Carlo (MC) methods can be used directly for this purpose, but they require a high computational effort \cite{Caflisch_1998, Niederreiter_1992}. Quasi Monte Carlo (QMC) and variance reduction techniques \cite{Caflisch_1998} may reduce the computational effort considerably without requiring much regularity. However, often we have high regularity in the stochastic variables, and this is not exploited by QMC methods.

Alternatives to Monte Carlo methods have been developed. For example, perturbation methods \cite{Kleiber_1992}, methods based on Neumann-series \cite{Babuska_2002}, or the spectral stochastic finite element method (SSFEM) \cite{Ghanem_1991, Ghanem_1999}. Stochastic Galerkin methods have been applied to various linear problems, see \cite{Ghanem_1991, Xiu_2002, Pulch_2014}. Nonlinear problems with stochastic loads have been tackled in \cite{Xiu_2002_2}. These Galerkin methods yield an explicit functional relationship between the independent random variables and the solution. In contrast with common MC methods, subsequent evaluations of functional statistics like the mean and covariance are very cheap. 

We consider an elliptic PDE in space including a random field as material parameters. The polynomial chaos approach and the stochastic Galerkin method yield a deterministic system of PDEs in space \cite{Pulch_2012}. In this paper, we introduce a spatial collocation technique  based on polynomial approximation by Lagrange interpolation. For the interpolation points we use a specific set of non-uniform points created by conformal maps, called Sinc points. Later, we use a small number of Sinc points as collocation points to compute a very accurate solution of the PDEs, see \cite{Maha_PhD}.

The paper is organized as follows: In Section \ref{sec:2}, we introduce a model problem, the structure of its polynomial chaos model and the stochastic Galerkin solution. In Section \ref{sec:3}, we illustrate the main theorem of Poly-Sinc approximation. In Section \ref{sec:4}, we review a Poly-Sinc collocation technique with the main collocation theorem. Finally, in Section \ref{sec:5}, we investigate numerical examples. We start with a simple example in one stochastic variable and then we discuss the general model from Section \ref{sec:2}.

\section{Stochastic Model Problem}
\label{sec:2}
In this paper, we are interested to solve the following stochastic partial differential equations:

\begin{equation}
\label{Eqn_01}
\begin{split}
\mathcal{L}(u)\equiv -\nabla \cdot (a(x,y,\Theta )\nabla u(x,y,\Theta))&=f(x,y) \text{ in } Q \times \Omega  \text{ and }\\
 u&=0 \text{ on } \partial Q \times \Omega,
\end{split}
\end{equation}

where \(\Theta =\left(\xi _1,\xi _2,\text{...},\xi _K\right)\) is a vector of stochastic parameters. These parameters are independent and uniformly distributed in \(I=[-1,1]\) and thus \(\Theta:\Omega\longrightarrow[-1,1]^K\) with an event space \(\Omega\). Moreover the domain of the spatial variables \(x\) and \(y\) is \(Q =(0,1)^2\). The function \(a(x,y,\Theta )\) is defined
as 

\begin{equation}
\label{Eqn_02}
a(x,y,\Theta )=a_0(x,y)+b_0\sum _{k=1}^K \xi _k a_k(x,y),
\end{equation}

where \(a_k\){'}s are functions in \(x\) and \(y\) only, \(b_0\) is a constant and, \(\xi _k\){'}s are the random variables. Without loss of generality, we consider \(a_0=1\) and \(b_0=1/2\). We assume that \(a(x,y,\Theta)\geq \alpha > 0\) for all \((x,y)\in Q\) and all \(\Theta \in [-1,1]^K\). Thus the differential operator in (\ref{Eqn_01}) is always uniformly elliptic.

In the rest of the section, we introduce the main concepts used in the solution of (\ref{Eqn_01}) with (\ref{Eqn_02}). Basically, we discuss the polynomial chaos in one- and multidimensional cases and the stochastic Galerkin method.

 
\subsection{Polynomial Chaos Expansion}
\label{PC}
Generalized Polynomial Chaos (gPC) is a particular set of polynomials in a given random variable, with which an approximation of a finite second-moment random variable is computed. This procedure is named Polynomial Chaos Expansion (PCE). This technique exploits orthogonal properties of polynomials involved, to detect a representation
of random variables as series of functionals.
Now, the function \(u\) can be expressed as an infinite series of orthogonal basis functions \(\Phi_{i}\) with suitable coefficient functions \(u_i\) as

\begin{equation}
\label{Eqn_030}
u(x,y,\Theta ) = \sum _{i=0}^{\infty} u_i(x,y)\Phi _{i}(\Theta ).
\end{equation}

The expansion in (\ref{Eqn_030}) converges in the mean square of the probability space. The truncation form including \(m+1\) basis functions leads to

\begin{equation}
\label{Eqn_03}
u(x,y,\Theta )\simeq {\widetilde{u}}(x,y,\Theta )= \sum _{i=0}^m u_i(x,y)\Phi _{i}(\Theta )
\end{equation}

\noindent with coefficients functions 

\begin{equation*}
\label{Eqn_031}
u_i(x,y)= \left\langle u(x,y,\Theta),\,\Phi _{i}(\Theta )\right\rangle, \,\,\, i=0,1,\ldots, m.
\end{equation*}
A fundamental property of the basis functions is the orthogonality,

\begin{equation}
\label{Eqn_07}
\left\langle \Phi_{i}(\Theta),\,\Phi_{j}(\Theta)\right\rangle=\int_{I^K}\Phi_{i}(\Theta)\,\Phi_{j}(\Theta)\,W(\Theta)d\Theta = c_{i}\,\delta_{i j}, \,\, \text{for all } {i}, {j},
\end{equation}
where \(c_{i}\) are real positive numbers and \(\delta_{i j}\) is the Kronecker-delta. In general, the inner product in (\ref{Eqn_07}) can be defined for different types of weighting function \(W\); however, it is possible to prove that the optimal convergence rate of a gPC model can be achieved when the weighting function \(W\) agrees to the joint probability density function (PDF) of the random variables considered in a standard form \cite{Xiu_2002,Witteveen_2006}. In this framework, an optimal convergence rate means that a small number of basis functions is sufficient to obtain an accurate PC model (\ref{Eqn_03}). Hence, the choice of the basis functions depends only on the probability distribution of the random variables $\Theta$, and it is not influenced by the type of system under study. In particular, if the random variables $\Theta$ are independent, their joint PDF corresponds to the product of the PDFs of each random variable: in this case, the corresponding basis functions $\Phi_{i}$ can be calculated as product combinations (tensor product) of the orthogonal polynomials corresponding to each individual random variable \cite{Eldred_2009, Pulch_2012, Xiu_2010}:

\begin{equation}
\label{Eqn_08}
\Phi_{i}(\Theta)=\Phi_{\mathbf{i}}(\Theta):= \prod^{K}_{r=1}\Phi_{i_r}^{(r)}(\xi_r),\,\, \mathbf{i}=(i_1,\ldots,i_K),
\end{equation}

\noindent where \(\Phi_{i_r}^{(r)}\) represents the univariate basis polynomial of degree \(i_r\) associated to the \(r\)th random parameter and with one-to-one correspondence between the integers \(i\) and the multi-indices \(\mathbf{i}\). We assume \({\rm{degree}}(\Phi_i) \leq {\rm{degree}}(\Phi_{i+1})\) for each \(i\). Now let

\begin{equation}
\label{Eqn_080}
\mathcal{R}_P=\left\{\Phi_{i}(\Theta): \sum^{K}_{r=1}i_r \leq P\right\},
\end{equation}

be the set of all multivariate polynomials up to total degree \(P\) as used in a Taylor expansion. Furthermore, for random variables with specific PDFs, the optimal basis functions are known and are formed by the polynomials of the Wiener-Askey scheme \cite{Xiu_2002}. For example, in the uniform probability distribution, the basis functions are the Legendre polynomials.

Using (\ref{Eqn_08}) and (\ref{Eqn_080}), it is possible to show that the total number of basis functions $m+1$ in (\ref{Eqn_03}) is expressed as
\begin{equation}
\label{Eqn_09}
m+1=\frac{(K+P)!}{K!P!}.
\end{equation}

\noindent The total degree of the PC (the maximum degree) $P$ can be chosen relatively small to achieve the desired accuracy in the solution.

In the case of the orthogonal polynomials, we can see that \(\Phi_0 (\Theta)=1\) and for orthonormal polynomials

\begin{equation}
\label{Eqn_10}
\left\langle \Phi_{i}(\Theta),\,\Phi_{i}(\Theta)\right\rangle=1.
\end{equation}
Once a PC model in the form of (\ref{Eqn_03}) is obtained, stochastic moments like the mean \(E(u)\) and the variance \(V(u)\) can be analytically calculated by the PC expansion coefficients as, see \cite{Xiu_2010},

\begin{eqnarray*}
\label{Eqn_11}
E(u(x,y,\Theta))& = &\int_{I^K} u(x,y,\Theta)\, W(\Theta)\, d\Theta \\
& = & \int_{I^K} u(x,y,\Theta)\,\Phi_{0}(\Theta) W(\Theta)\, d\Theta \\
& = & \left\langle u(x,y,\Theta), \Phi_{0}(\Theta)\right\rangle.
\end{eqnarray*}

Using the PC expansion of $u(x,y,\Theta)$ given in (\ref{Eqn_030}) and the orthogonality of the basis functions $\Phi_{i}(\Theta)$ to get

\begin{equation*}
\label{Eqn_12}
E\left(u(x,y,\Theta)\right)=u_0(x,y).
\end{equation*}

The variance can be derived by

\begin{align*}
\label{Eqn_13}
V(u(x,y,\Theta))& = \int_{I^K} [u(x,y,\Theta)-E(u(x,y,\Theta))]^2\, W(\Theta)\, d\Theta \\
& = \resizebox{.8\hsize}{!}{$\int_{I^K} [u^2(x,y,\Theta)+E^2(u(x,y,\Theta))-2u(x,y,\Theta)E(u(x,y,\Theta))]\, W(\Theta)\, d\Theta$} \\
& = \int_{I^K} [u^2(x,y,\Theta)+u^2_0(x,y)-2 u_0(x,y)\,u(x,y,\Theta)]\, W(\Theta)\, d\Theta\\
& = \left\langle u(x,y,\Theta), u(x,y,\Theta)\right\rangle - u^2_0(x,y).
\end{align*}

Again, use the PC expansion in (\ref{Eqn_03}) and orthonormal polynomials basis satisfying (\ref{Eqn_10}), to get

\begin{equation*}
\label{Eqn_14}
V(u(x,y,\Theta)) \approx \sum^{m}_{i=1}u^{2}_{i}(x,y).
\end{equation*}

It is clear now that, in order to obtain a PC model in (\ref{Eqn_03}) and the stochastic moments, the coefficients functions $u_i (x,y)$ must be computed. The PC coefficient estimation depends on the type of the resulting system from the chaos expansion, not only the PC truncation. 
 
\subsection{Stochastic Galerkin Method}

To solve the problem in (\ref{Eqn_01}) and (\ref{Eqn_02}), a Galerkin method is used along side the PC. The main idea is to assume that the solution of (\ref{Eqn_01}) and (\ref{Eqn_02}) is written as expansion in (\ref{Eqn_03}) and then use the PC theory introduced in the previous section. This process transform the SPDE (\ref{Eqn_01}) and (\ref{Eqn_02}) into a deterministic system of PDEs. 

To recover the coefficient functions \(u_i(x,y)\) we apply the inner product of (\ref{Eqn_01}) with the basis polynomial \(\Phi _{j}(\Theta )\)

\begin{equation}
\label{Eqn_04}
\left\langle  \mathcal{L}(\widetilde{u}) - f(x,y),\Phi _{j}(\Theta )\right\rangle =0 \,\,\,\,\, \text{for } j=0,1,\ldots,m.
\end{equation}

Substituting (\ref{Eqn_03}) in (\ref{Eqn_01}) we obtain 

\begin{equation*}
\label{Eqn_05}
\mathcal{L}(\widetilde{u})=-\nabla \cdot \left(\nabla \sum _{i=0}^m u_i(x,y)\Phi _{i}(\Theta )\right)-\frac{1}{2}\sum _{k=1}^K \xi _k \nabla \cdot \left(a_k\nabla \sum _{i=0}^m u_i(x,y)\Phi _{i}(\Theta )\right).
\end{equation*}

Now applying the inner residual product in (\ref{Eqn_04}) and use the orthogonality property of the multivariate basis \(\Phi _{i}\){'}s to get

\begin{equation}
\label{Eqn_06}
-\nabla ^2u_j-\frac{1}{2}\sum _{k=1}^K \sum _{i=0}^m \left\langle  \xi _k \Phi _{i}(\Theta ),\Phi _{j}(\Theta )\right\rangle  \nabla \cdot \left(a_k\nabla u_i\right)=F_j
\end{equation}

where \(F_j(x,y) = \left\langle  f(x,y),\Phi _{j}(\Theta )\right\rangle\) forms an \((m+1)\) vector and the array \(\left\langle  \xi _k\Phi _{i},\Phi _{j} \right\rangle\) is a triple tensor of dimension \(K \times (m+1)\times (m+1)\). (\ref{Eqn_06}) is a system of elliptic PDEs with unknown variables \(u_i(x,y)\), \(i=0,1,2,\ldots,m\). With large number of random variables \(K\) (say \(K>4\)) the size of the system in (\ref{Eqn_06}) becomes huge due to (\ref{Eqn_09}). One of our targets in the solution of the system in (\ref{Eqn_06}) is to use a collocation method to achieve a high accuracy with small numbers of collocation points. The proposed method in this report is to use Sinc points in a Lagrange interpolation. 


\subsection{Quadrature}
The inner product \(\left\langle ., . \right\rangle\) is defined by an integral. For the integration of polynomials analytic methods are used. Alternatively, we can use highly accurate quadrature techniques  to evaluate the integrals exactly except for round-off errors. We omit the details of these techniques, since they can be easily found in several textbooks. For example, descriptions of Gaussian quadrature can be found in most texts on numerical analysis \cite{Stoer_2002}, while \cite{Stenger_2010} contains descriptions of Sinc quadratures over finite, semi-infinite, infinite intervals and contours.

\section{Poly-Sinc Approximation}
\label{sec:3}
In this section we introduce the Lagrange approximation at Sinc points as interpolation points. This approximation is called Poly-Sinc approximation \cite{Maha_PhD}. It was first introduced to provide a uniform approximation for a function and its derivatives as well \cite{Stenger_2013}. In \cite{Maha_PhD}, the main results of this approximation have been extended and have been used to solve differential equations.


Given a set of data \( \left\{x_k, u(x_k)\right\}_{k=-M}^N\) where the \(x_k\) are interpolation points on \((a,b)\). Then there is a unique polynomial \(P_n(x),\, n=M+N+1\) of degree at most \(n-1\) satisfying the interpolation condition,

$$
P_n(x_k)=u_k,\, \, k=-M,\, ... ,N.
$$

\noindent In this case \(P_n(x)\) can be expressed by the Lagrange polynomials as

\begin{equation*}
\label{eq:11}
P_n(x)=\sum _{k=-M}^N b_k(x) \, u(x_k),
\end{equation*}

\noindent with,

\begin{equation*}
\label{eq:12}
b_k(x)=\frac{g(x)}{(x-x_k)\,g'(x_k)},\,\,\, g(x)=\prod_{j=-M}^{N} \left(x-x_j\right).
\end{equation*}

Now \(x_k\)'s are Sinc points on \((a,b)\) defined as \cite{Stenger_2010}

\begin{equation}
\label{eq:sp}
 x_k=\frac{a+b\, {\rm{e}}^{kh}}{1+{\rm{e}}^{kh}}.
\end{equation}

Corresponding to such a scheme, we define a row vector \(\bm{B}\) of basis functions and an operator \(V_mu\) that maps a function \(u(x)\) into a column vector of dimension \(n=M+N+1\) by

\begin{align*}
\label{equation:13}
\bm{B}(x) &=\left(b_{-M,h}(x),\, \ldots ,\, b_{N,h}(x)\right)\\
V_n u &=\left(u\left(x_{-M}\right),\, \ldots ,\, u\left(x_N\right)\right)^{\top}.
\end{align*}

\noindent This notation enables us to write the above interpolation scheme in simple operator form, as

\begin{equation}
u(x)\simeq \bm{B}(x)  \, V_n u.
 \label{equation:15}
\end{equation}

This approximation, like regular Sinc approximation, yields an exceptional accuracy in approximating the function that is known at Sinc points, \cite{Stenger_2010}. Unlike Sinc approximation, it gives a uniform exponential convergence rate when differentiating the interpolation formula given in (\ref{equation:15}), see \cite{Stenger_2013}. Next, we assume that $M=N$ and that $n=2N+1$ is the total number of Sinc points. Then the upper bound of error for Poly-Sinc approximation is given as

\begin{equation}
\label{eq:13}
\underset{x\in (a,b)} {\sup}|u(x)-\bm{B}(x) \, V_n u|\leq A\frac{\sqrt{N}}{B^{2N}}\, \, \exp  \left(\frac{-\pi^2  N^{\frac{1}{2}}}{2}\right),
\end{equation}

\noindent where \(A>0\) and \(B>1\) are two constants, independent of \(N\). For the proof of (\ref{eq:13}), see \cite{Stenger_2013}. 


Another criterion to discuss the convergence and stability of the Poly-Sinc approximation is the Lebesgue constant. In \cite{Youssef_2016} an estimate for the Lebesgue constant for Lagrange approximation at Sinc points (\ref{eq:sp}) has been derived as

\begin{equation*}
\label{eq:6}
\Lambda _{n} \approx \frac{1}{\pi }\log (n+1)+1.07618.
\end{equation*}

Next we extend these results from the one-dimensional case to the multi-dimensional case. 

Let $X=(x_1,.....,x_l)$ be a point in $Q=[a,b]^l$, then Lagrange approximation of a function $u(X)$ can be defined by a nested operator as

\begin{equation}
\label{eq:800}
       (P_n u)(X)=\sum_{k_1=-M_1}^{N_1} \sum_{k_2=-M_2}^{N_2} \ldots \sum_{k_l=-M_l}^{N_l} u(X_{\pmb{k}})\,  b_{k_1}(x_1)b_{k_2}(x_2) \ldots b_{k_l}(x_l), 
\end{equation}
where \(u(X_{\pmb{k}})=U=u(x_{1,k_1}, \ldots,x_{l,k_l} )\) with $k_i=-M_i,\ldots,N_i$. We can write the approximation (\ref{eq:800}) in the operator form

\begin{equation}
u(X)\simeq \bigodot^{l}_{i=1} {\bm{B}}_i(X) \, U,\,\,\, i=1,2,\ldots , l.
 \label{eq:155}
\end{equation}

Next, we assume \(M_i=N_j=N,\, i,j=1,\ldots,l\) and \(n=2N+1\) is the number of Sinc points in each dimension \(i=1,2,\ldots,l\). The convergence and stability of the approximation (\ref{eq:155}) are discussed in \cite{Youssef_2016_3} and \cite{Youssef_2016}. For the upper bound of the error \(E_n\), we have

\begin{equation}
\label{eq:8011}
E_n = \underset{X\in Q} {\sup} \, |u(X)- \bigodot^{l}_{i=1} {\bm{B}}_i(X) \, U|\leq \sum_{i=0}^{l-1} {\left(C_i \log^i{N}\right)\frac{\sqrt{N}}{\gamma_i^{2N}}\, \exp  \left(\frac{-\pi^2  N^{\frac{1}{2}}}{2}\right)},
\end{equation}

\noindent where \(C_i>0\), \(\gamma_i>1\), \(i=1, \ldots, l\) are two sets of constants, independent of \(N\).

The notation $\Lambda_{n,l}$ is used to denote the Lebesgue constant using $n$ interpolation points in each dimension $i=1,2,\ldots,l$, i.e. $n^l$  Sinc points in total. If \(P_n(X)\) is defined as in (\ref{eq:800}), then:

\begin{equation}
\label{eq:805}
\Lambda _{n,l}\leq \left(\frac{1}{\pi }\log (n+1)+1.07618\right)^l .
\end{equation}

\section{Poly-Sinc Collocation Method}
\label{sec:4}
In \cite{Youssef_2016_1}, a collocation method based on the use of bivariate Poly-Sinc interpolation defined in (\ref{eq:155}) is introduced to solve elliptic equations defined on rectangular domains. In \cite{Youssef_2016_4}, Poly-Sinc collocation domain decomposition method for elliptic boundary value problems is investigated on complicated domains. The idea of the collocation method is to reduce the boundary value problem to a system of algebraic equations which have to be solved subsequently. To start let us introduce the following collocation theorem.

\begin{thm}
Let \(u:\overline{Q}\rightarrow \mathbb{R}\) be an analytic bounded function on the compact domain $\overline{Q}$. Let $U=\left\{u(x_j,y_k)\right\}^{N}_{j,k=-N}$ be a vector, where $x_j$ and $y_k$ are the Sinc points. If ${\widetilde{U}}=\left\{{\widetilde{u}_{j k}}\right\}^{N}_{j,k=-N}$ is a vector satisfying
 
$$\left\|U-\widetilde{U}\right\|_{\infty}=\max_{j,\,k}\,\left|u_{j k}-\widetilde{u}_{j k}\right| < \delta,$$
then
\begin{equation}
\label{eq:coler}
\left\|u(x,y)- \bigodot^{2}_{i=1} {\bm{B}}_i(x,y)\, \widetilde{U}\right\| < E_n \, + \, \delta\, \Lambda_{n,2}, 
\end{equation}
where \(n=2N+1\), $E_n$ from (\ref{eq:8011}), and $\Lambda_{n,2}$ from (\ref{eq:805}). 

\end{thm}

\newpage

\begin{prf}
We apply triangle inequality 

\begin{align*}
\left\| u(x,y)-\bigodot^{2}_{i=1} {\bm{B}}_i(x,y)\, \widetilde{U}\right\| & \leq  
\left\|u(x,y)-\bigodot^{2}_{i=1} {\bm{B}}_i(x,y)\, U \right\| \\
  & \,\,\, + \left\|\bigodot^{2}_{i=1} {\bm{B}}_i(x,y)\, U-\bigodot^{2}_{i=1} {\bm{B}}_i(x,y)\, \widetilde{U} \right\| \\
& \leq  E_n + \delta \, \left\|\bigodot^{2}_{i=1} {\bm{B}}_i(x,y)\right\|\\
&\leq  E_n + \delta \, \Lambda_{n,2}, 
\end{align*}
which is the statement of the theorem.
$\hfill\square$
\end{prf}

This theorem guarantees an accurate final approximation of \(u\) on its domain of definition provided that we know a good approximation to \(u\) at the Sinc points.

To set up the collocation scheme, let us consider the following partial differential operator,

\begin{align}
\label{eq:1}
{\cal{L}}u \equiv u_{x\,x}+u_{y\,y} & =f(x,y),\,\, (x,y)\in Q,\\
 u(x,y) & =u_{ex}(x,y),\,\, (x,y)\in \partial Q, \nonumber
\end{align}

\noindent where $Q=\left\{a<x<b,\, c<y<d \right\}$ and $u_{x\,x}=\frac{\partial^2 u}{\partial{x^2}}$, $u_{y\,y}=\frac{\partial^2 u}{\partial{y^2}}$.

The first step in the collocation algorithm is to replace \(u(x,y)\) in Eq. (\ref{eq:1}) by the Poly-Sinc approximation  defined in (\ref{eq:155}). Next, we collocate the equation by replacing $x$ and $y$ by Sinc points

$$x_i=\frac{a+b\, {\rm{e}}^{i\, h}}{1+ {\rm{e}}^{i\, h}},\, i=-M,\ldots,N$$ and 
$$y_q=\frac{c+d \, {\rm{e}}^{q\, h}}{1+{\rm{e}}^{q h}},\, q=-M,\ldots,N.$$

In this case, we have,
 
 \begin{equation*}
u_{x\,x}(x_i,y_q)\approx \sum_{k=-M}^{N}\sum _{j=-M}^N u_{j k}\, B^{''}(j,h)(x_i) B(k,h)(y_q),
\end{equation*}

where,
\begin{equation*}
\label{eq:2000}
B(j,h)(x_i)=\delta_{j\,i}=\begin{cases}
  0 &   j\neq i. \\
  1 &   j=i,
\end{cases}
\end{equation*}
 
 and $B^{''}(j,h)(x_i)$ defines an $n \times n$ matrix, with \(n=M+N+1\)

\begin{equation*}
\label{eq:2001}
B^{''}(j,h)(x_i)=[b_{j i}]=\begin{cases}
 \frac{-2 g'(x_i)}{(x_i-x_j){}^2g'(x_j)}+\frac{g''(x_i)}{(x_i-x_j)g'(x_j)}
& \text{if}\, j\neq i \\
\\
\mathlarger{\sum}^N_{s=-M} \, \mathlarger{\sum}^N_{\substack{l=-M\\
l,s\neq i}} \frac{1}{(x_i-x_l)(x_i-x_s)} & \text{if}\, j=i.
\end{cases}
\end{equation*}

So, 

\begin{equation*}
\label{eq:2002}
{\cal{U}}_{x\,x} = \left(u_{x\,x}(x_i,y_q)\right)_{i,q=-M,\ldots,N}= {\cal{M}}_{1} \, {\cal{U}},
\end{equation*}
 
\noindent where ${\cal{M}}_{1}$ is a $n^2 \times n^2$ matrix defined as,
\begin{equation*}
\label{eq:2003}
{\cal{M}}_{1}=\begin{cases}
b_{j \, i} & k=q \wedge \, i,\, j,\, k,\, q=-M,\,...,\,N\\
\\
0 & k\neq q \wedge \, i,\, j,\, k,\, q=-M,\,...,\,N,
\end{cases}
\end{equation*}

\noindent and where ${\cal{U}}_{x\,x}$ is collected in a vector of of length $n^2$. Likewise, it holds that

\begin{equation*}
\label{eq:2004}
{\cal{U}}_{y\,y} = \left(u_{y\,y}(x_i,y_q)\right)_{i,q=-M,\ldots,N}= {\cal{M}}_{2} \, {\cal{U}},
\end{equation*}

where ${\cal{M}}_{2}$ is defined in the same way as ${\cal{M}}_{1}$.

 The differential equation has been transformed to a system of $n^2$ algebraic equations,
\begin{equation*}
 {\cal{A}} \, \mathcal{U}=\mathcal{F}, 
\label{eq:fin}
\end{equation*}

\noindent where $\mathcal{U}$ is the vector of length $n^2$ including the unknowns $u_{i\, q}$ and
 
 \begin{equation*}
 \label{eq:2005}
 {\cal{A}}={\cal{M}}_{1} +{\cal{M}}_{2}.
 \end{equation*}
 
 The right hand side $\mathcal{F}$ is a vector of Length $n^2$ and defined as
 
 \begin{equation*}
 \label{eq:2006}
 \mathcal{F}=f(x_i,\,y_q), \, i,\, q=-M, ..., N. 
 \end{equation*}

Now the PDE (\ref{eq:1}) is transformed to a system of $n^2$ algebraic equations in $n^2$ unknowns. The boundary conditions are collocated separately to yield $4 n $ algebraic equations. More precisely, 

\begin{align*}
u(a,y_j) & =u_{ex}(a,y_j) \\
u(b,y_j) & =u_{ex}(b,y_j) \\
u(x_i,c) & =u_{ex}(x_i,c) \\
u(x_i,d) & =u_{ex}(x_i,d),
\end{align*}

where $x_i$ and $y_j$ are the Sinc points defined on $(a,b)$ and $(c,d)$, respectively. Adding these $4n$ equations to the $n^2 \times n^2$ algebraic system, produced from the collocation of the PDE, yields a rectangular system of linear equations. Finally, solving this least squares problem yields the desired numerical solution.

\newpage
Note:
\begin{itemize}
\item In our calculations we used a multiplier factor $\tau=10^3$ in the collocation steps of the homogenous boundary conditions. This factor emphasizes the boundary values and improve the error behavior at the boundaries.

\item The Poly-Sinc collocation technique is based on the collocation of the spatial variables using Sinc points. This means that it is valid also for PDEs with space-dependent coefficients. Moreover, it can be generalized to solve a system of PDEs.

\end{itemize}


\section{Numerical Results}
\label{sec:5}
In this section, we present the computational results. Mainly, we discuss two examples. The first simple example includes one stochastic parameter. In the second example we solve the model problem introduced in Section \ref{sec:2}.


\subsection{One Stochastic Variable}
Consider the Poisson equation in two spatial dimensions with one random parameter. This problem is described by the following SPDE

\begin{eqnarray}
\label{Eqn_100}
a(\xi)\left(u_{xx}(x,y,\xi)+u_{yy}(x,y,\xi)\right) & = & f(x,y) \,\, \text { on } Q\times \Omega \\
u(x,y,\xi) & = & 0 \,\, \text { on } \partial Q\times \Omega, \nonumber
\end{eqnarray}

where $Q=(-1,1)^2$ is the spatial domain and $\Omega$ is an event space and $\xi:\Omega \rightarrow [-1,1]$ is a random variable. The function \(a(\xi)=\xi + 2\) is a linear function of a uniformly distributed random variable $\xi$ and \(f(x,y)=1\) for all $(x,y)\in Q$.

Now, we use the PC representation in (\ref{Eqn_03}) with $m=3$ to have

\begin{equation}
\label{Eqn_101}
u(x,y,\xi)=\sum^{3}_{i=0}{u_i(x,y)\,\Phi_i(\xi)},
\end{equation}

\noindent where $\Phi_i$'s are the univariate orthonormal Legendre polynomials defined on $[-1,1]$. Substitution of (\ref{Eqn_101}) in the SPDE (\ref{Eqn_100}) yields the residual

\begin{equation*}
\label{Eqn_102}
R=(\xi + 2) \sum^{3}_{i=0}{\left((u_i)_{xx}+(u_i)_{yy}\right)\,\Phi_i(\xi)} - 1.
\end{equation*}

We then perform a Galerkin projection and use the orthogonality of Legendre polynomials, which yields the system of elliptic PDEs

\begin{align}
\label{Eqn_103}
\sum^{3}_{i=0}{\left\langle \Phi_k, \, (\xi + 2)\, \Phi_i \right\rangle}\mathcal{L}{u_i} & =\left\langle 1, \Phi_k\right\rangle  \text{   for } k=0,1,2,3,\,\, \text{on } Q\\
u_i & =  0  \text{ for } i=0,1,2,3\,\, \text{on } \partial Q, \nonumber
\end{align}

\noindent where $\mathcal{L}{u_i}=(u_i)_{xx}+(u_i)_{yy}$. It holds that $\left\langle 1, \Phi_k \right\rangle=\delta_{1 k}$.

\noindent The computational results of this example are given in the following experiments.

\begin{figure}[H]
\centering
\begin{subfigure}{.5\textwidth}
  \centering
  \includegraphics[width= \linewidth]{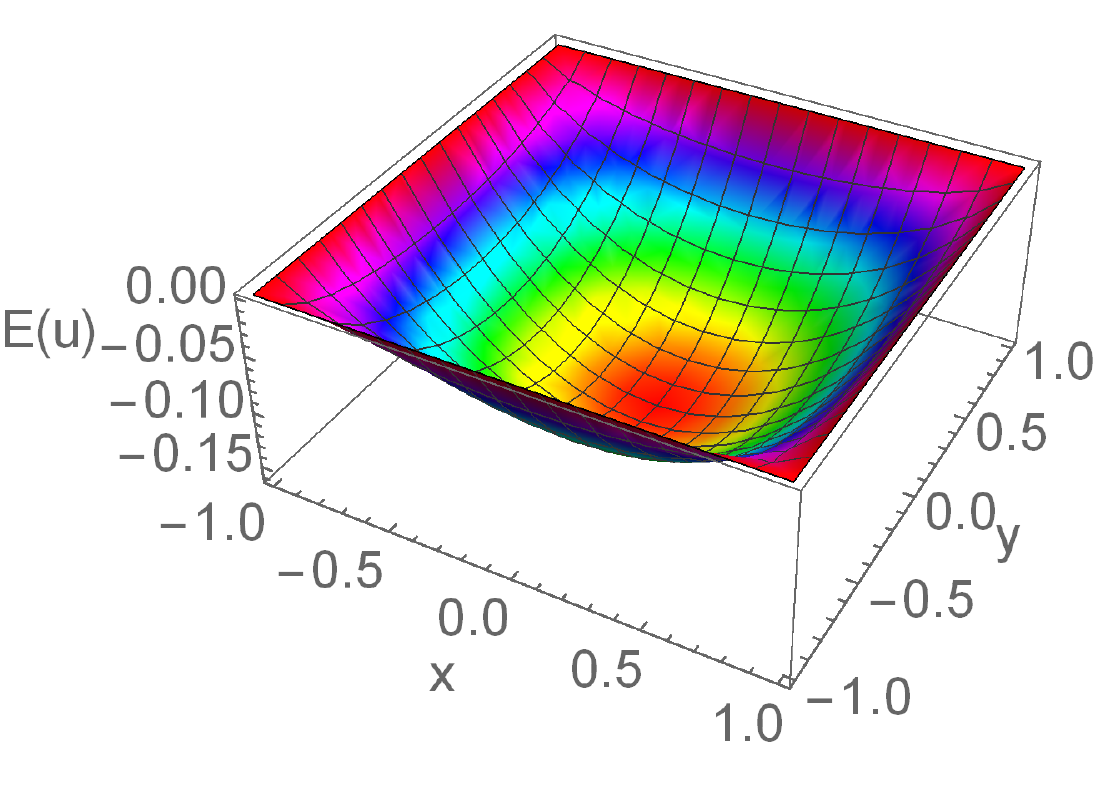}
  \caption{$E(u)$.}
\end{subfigure}%
\begin{subfigure}{.5\textwidth}
  \centering
  \includegraphics[scale=0.35]{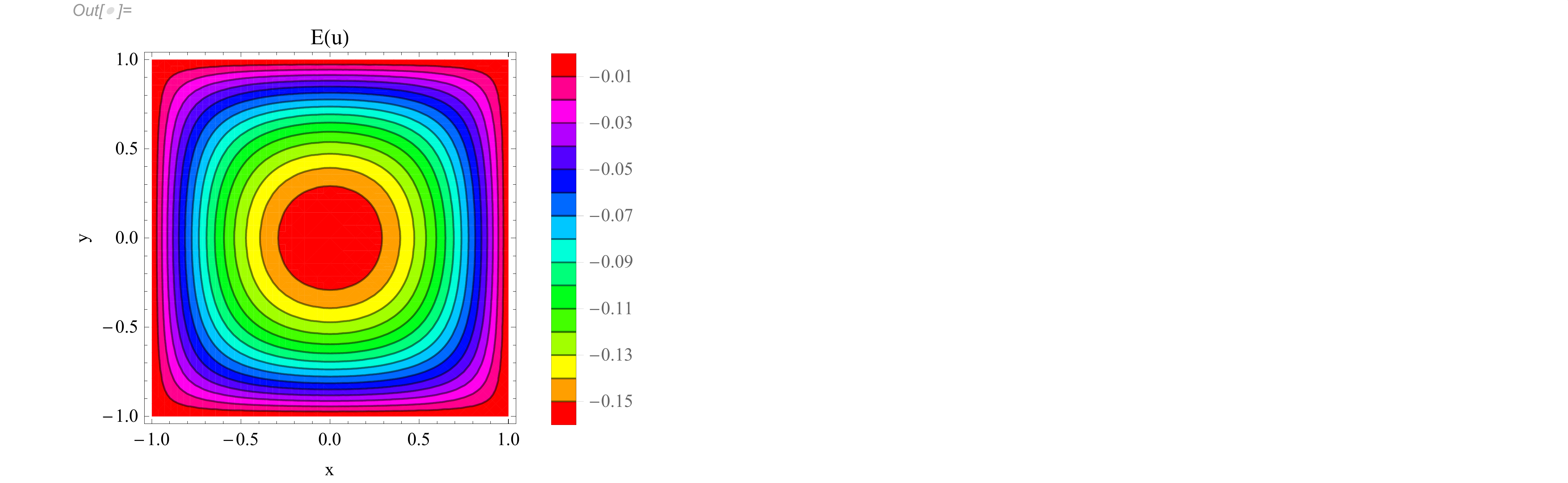}
  \caption{Contour plot of $E(u)$}
  \label{fig:sub2}
\end{subfigure}
\caption{The expectation, $E(u)$, using $m=3$ and Poly-Sinc with $N=5$.}
\label{fig:01}
\end{figure}

\begin{figure}[H]
\centering
\begin{subfigure}{.5\textwidth}
  \centering
  \includegraphics[width= \linewidth]{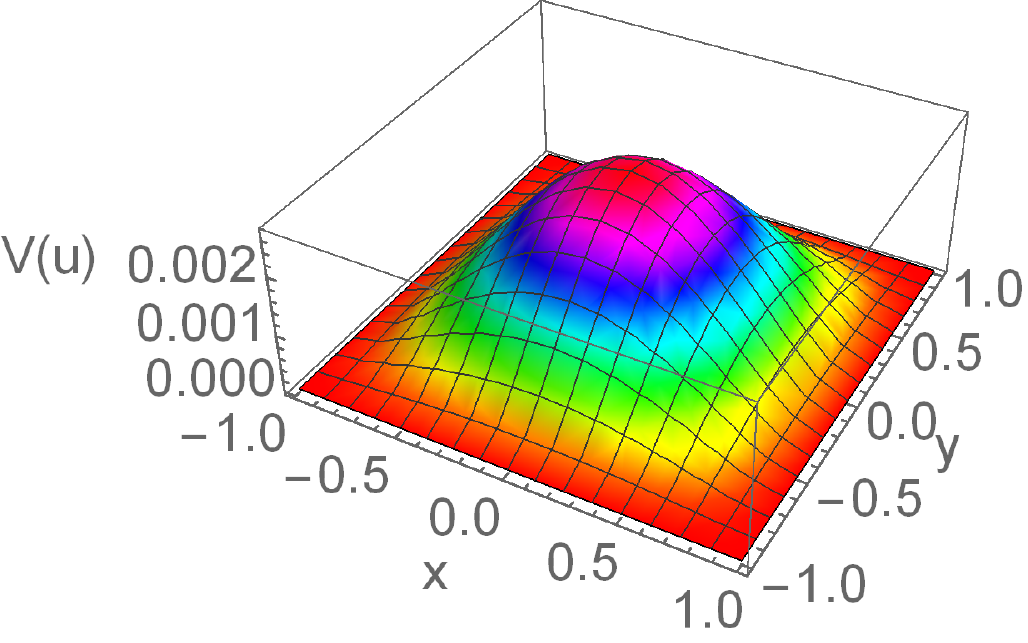}
  \caption{$V(u)$.}
\end{subfigure}%
\begin{subfigure}{.5\textwidth}
  \centering
  \includegraphics[width=.9\linewidth]{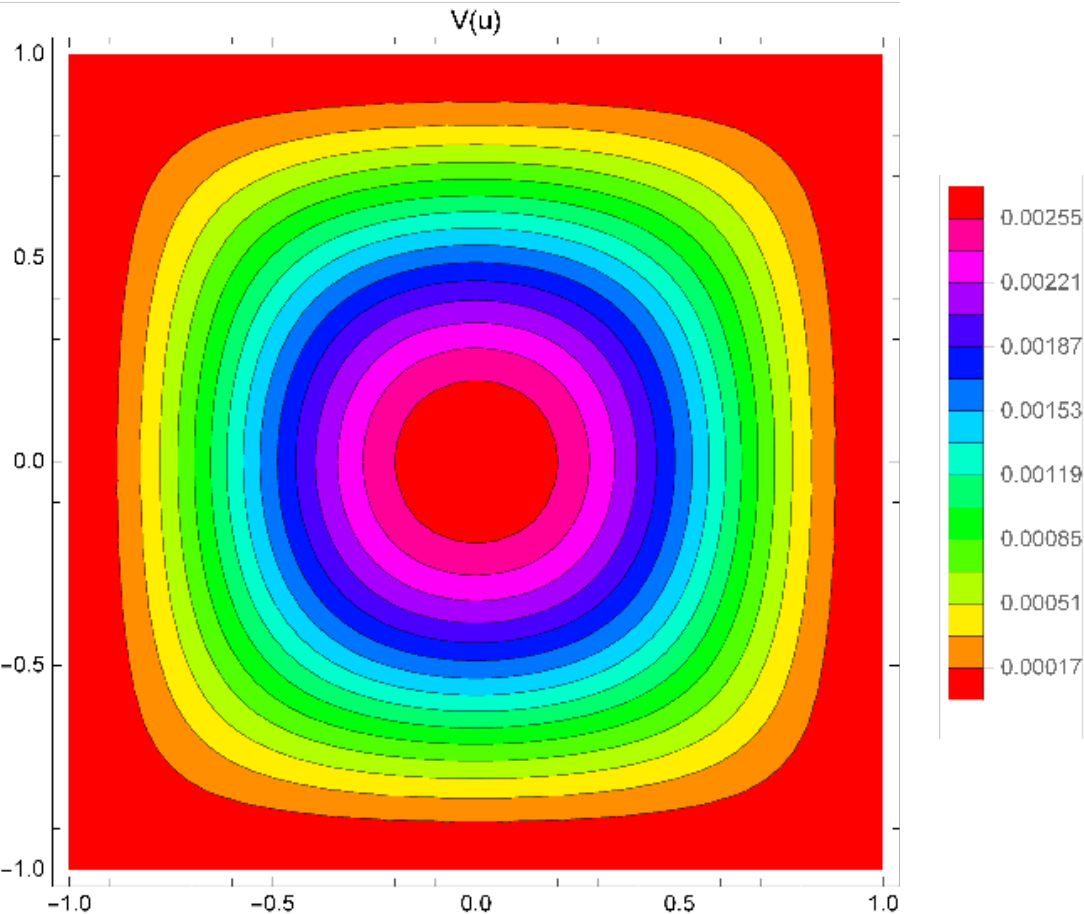}
  \caption{Contour plot of $V(u)$}
  \label{fig:sub2}
\end{subfigure}
\caption{The variance, $V(u)$, using $m=3$ and Poly-Sinc with $N=5$.}
\label{fig:02}
\end{figure}

\begin{expt}{$\mathbf{E(u)}$ \bf{and} $\mathbf{V(u)}$}
\newline
In this experiment, we use Poly-Sinc collocation from Section \ref{sec:4} to solve the system of PDEs in (\ref{Eqn_103}). In our computation, we use $N=5$, i.e. $11 \times 11$ of 2D grid of Sinc points defined as in (\ref{eq:sp}) on the domain $Q$. As a result of the Poly-Sinc solution, the coefficient functions $u_i(x,y)$ are obtained. In Fig. \ref{fig:01}, the expectation $E(u)=u_0(x,y)$ and its contour plot are represented while in Fig.~\ref{fig:02}, the variance calculations are presented. 
\end{expt}

\begin{figure}[H]
\centering
\includegraphics[scale=0.8]{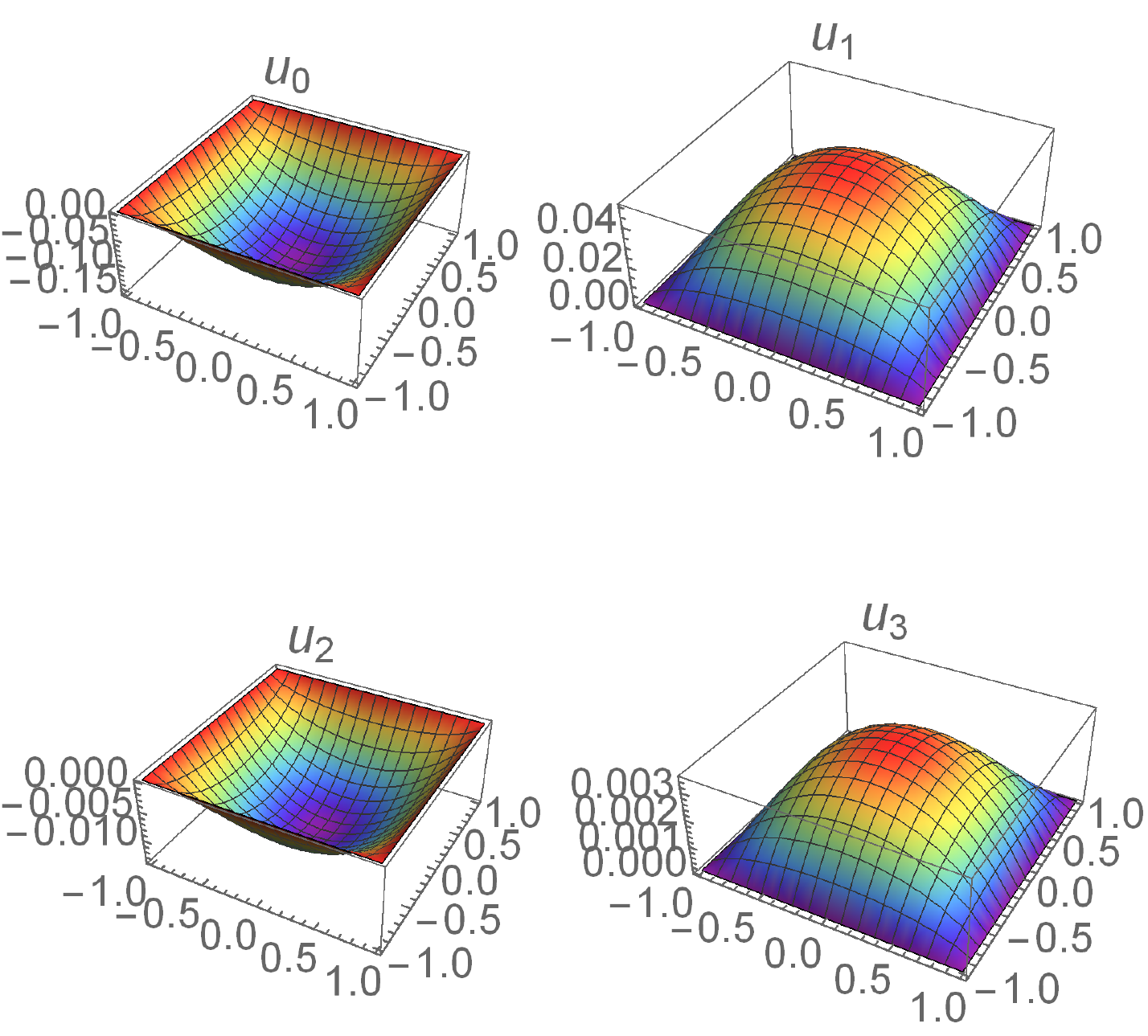}
\caption{Coefficients functions, $u_i (x,y),\, i=0,\ldots,3$.}
\label{fig:func1}
\end{figure}

\begin{figure}[H]
\centering
\includegraphics[scale=0.7]{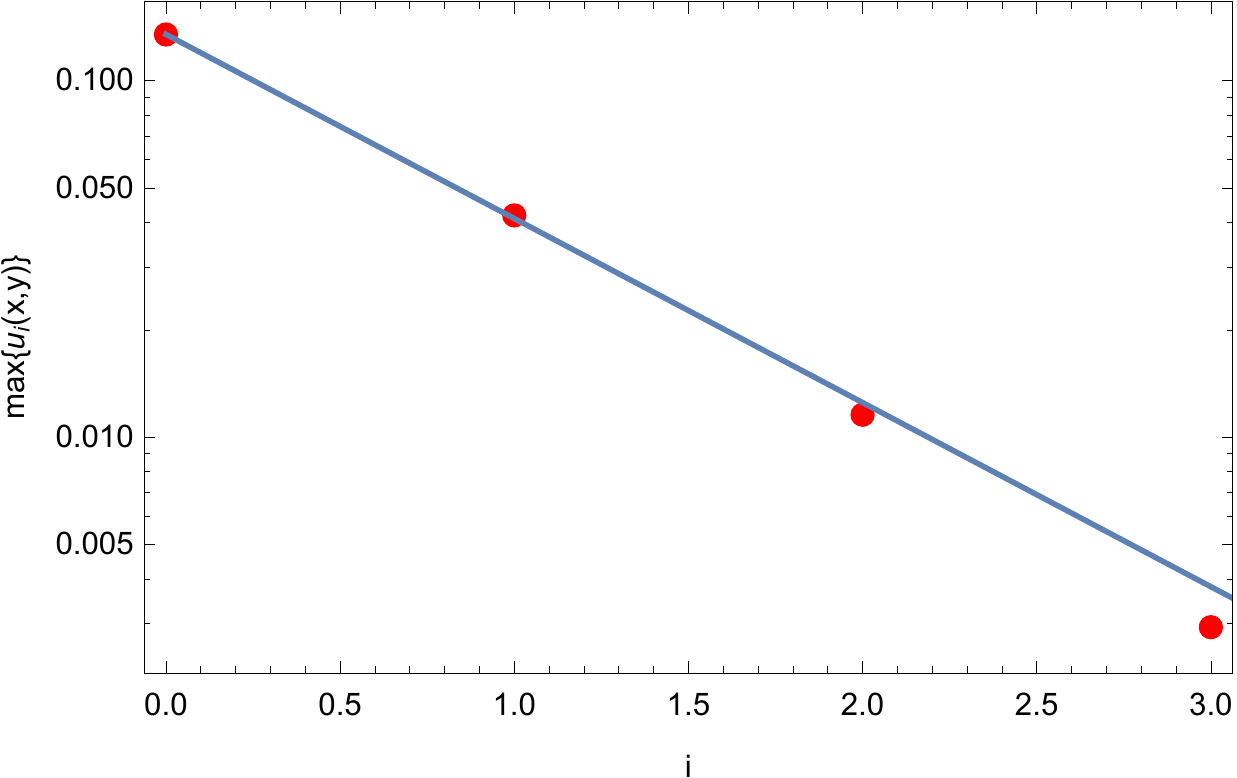}
\caption{Logarithmic plot of maximum of coefficient functions $u_i,\,\, i=0,\ldots,3$. The dots are the calculated maximum and the solid line represent the exponential fitting function $0.135\,e^{-1.2 \, i}$.}
\label{fig:03}
\end{figure}

\newpage
\begin{expt}{\bf{Coefficients functions}}
\newline
As we mentioned above, to get an accurate result, just a small number of orthogonal polynomials, $\Phi_i$, is needed. In our computations, we used $m=3$, i.e. four orthonormal Legendre polynomials. The $4$ coefficients functions, $u_i (x,y),\, i=0,\ldots,3$, are given in Fig.\ref{fig:func1}. In addition, we verify that this number is sufficient by showing that the coefficient functions $u_i$ tend to zero as $m$ increases. The results are given in Fig.\ref{fig:03}. In Fig.\ref{fig:03}, the dots represent the maximum of the coefficient functions $u_i(x,y)$ on the spatial domain. We then use these maximum values in a least square
estimation to find the coefficients of the decaying rate function $\alpha \, \exp (-\beta s)$, where $\alpha$ and $\beta$ are constants. In Fig.\ref{fig:03}, the solid line represents the best fitting function with $\alpha=0.14$ and $\beta=1.2$. This means that the coefficient functions $u_i(x,y)$ follow an exponentially decay relation.
\end {expt}

\begin{figure}[H]
\centering
\begin{subfigure}{.5\textwidth}
  \centering
  \includegraphics[width= \linewidth]{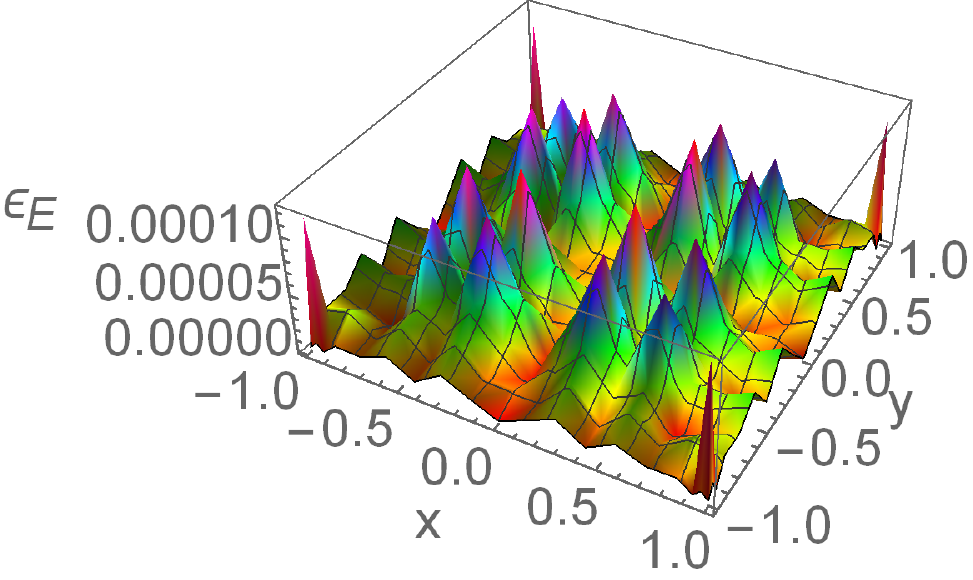}
  \caption{Absolute error in $E(u)$.}
\end{subfigure}%
\begin{subfigure}{.5\textwidth}
  \centering
  \includegraphics[width=.9\linewidth]{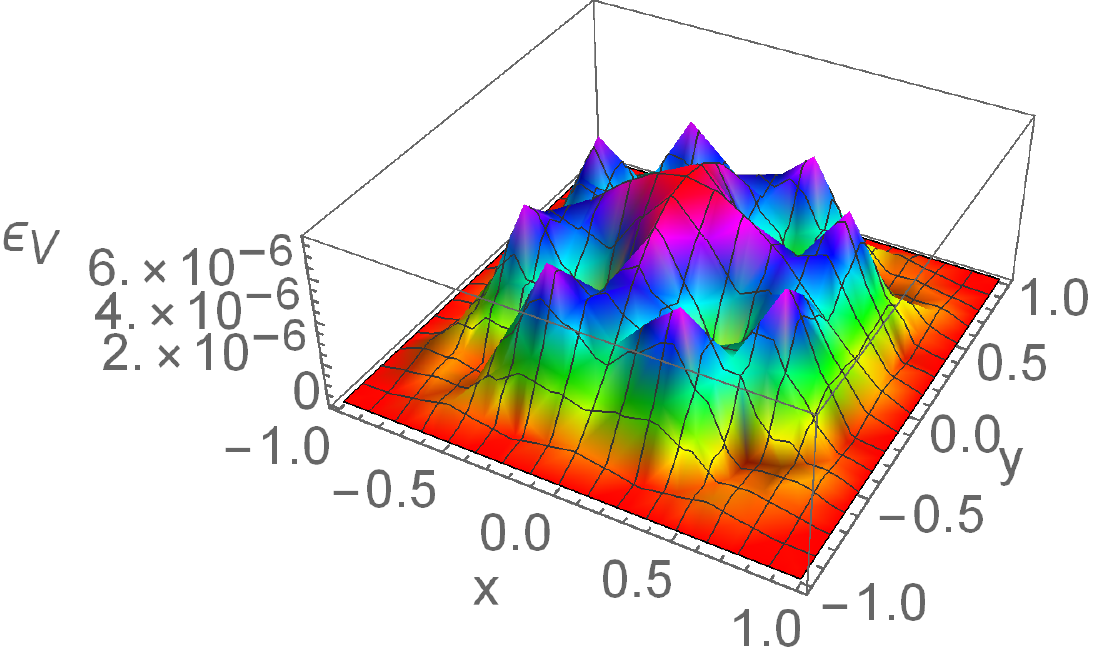}
  \caption{Absolute error in $V(u)$}
  \label{fig:sub2}
\end{subfigure}
\caption{Absolute error between the Poly-Sinc calculation and the calculations obtained from $100$ solutions.}
\label{fig:04}
\end{figure}

\begin{figure}[H]
\centering
\includegraphics[width= \linewidth]{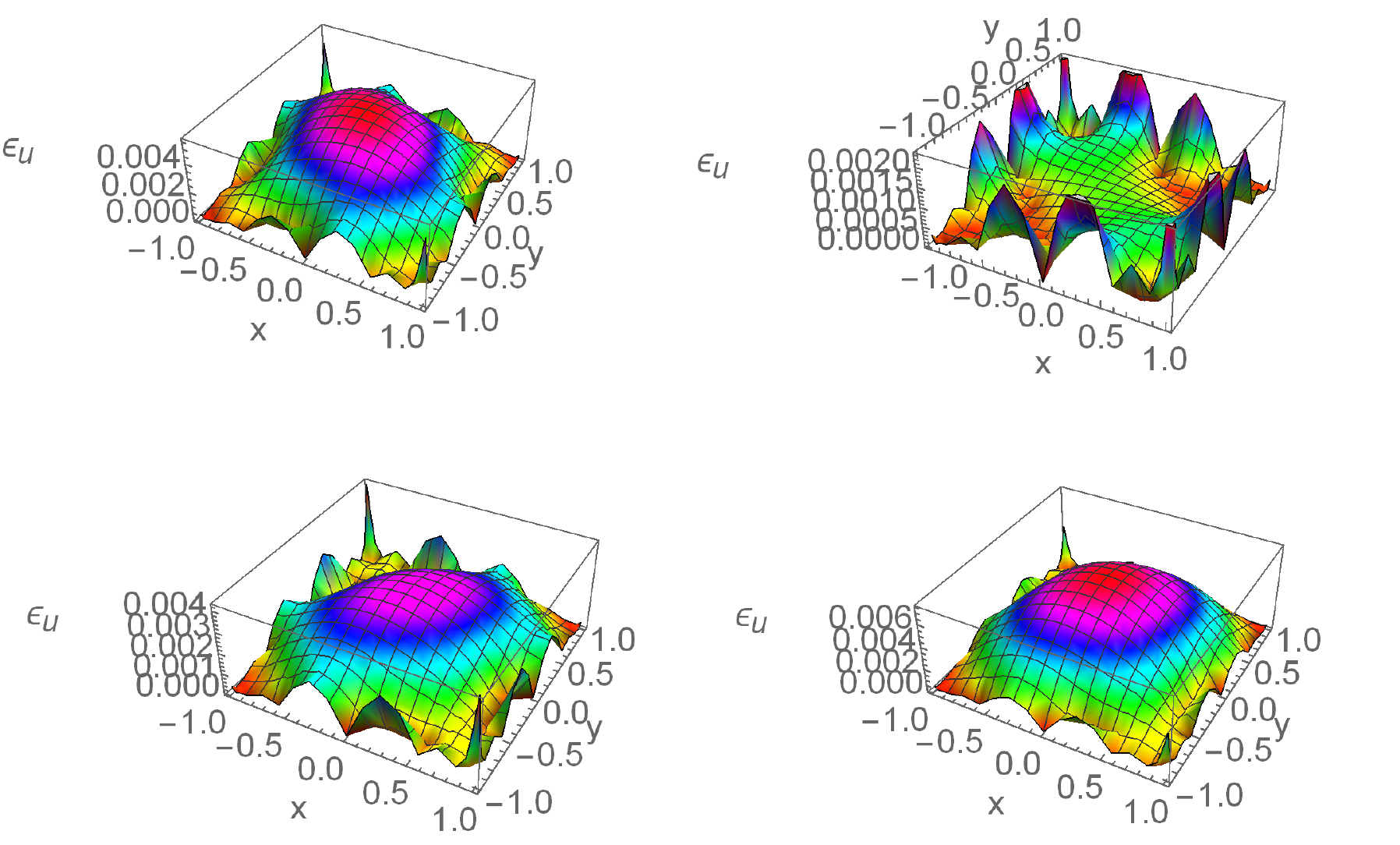}
\caption{Absolute error in $u$ for some discrete $\xi \in \left\{-0.757, 0, 0.757, 0.989\right\}$.}
\label{fig:05}
\end{figure}

\begin{expt}{\bf{Error}}
\newline
To discuss the convergence of Poly-Sinc solution, we need a reference (nearly exact) solution. For that, we create a discrete list of PDEs of the equation (\ref{Eqn_100}) at a finite set of instances of $\xi \in [-1,1]$. We choose $100$ points of Gauss-Legendre nodes as values of $\xi \in [-1,1]$ and create corresponding $100$ PDEs. To solve each one of these $100$ equations we use Mathematica Package NDSolve. NDSolve uses a combination of highly accurate numeric schemes to solve initial and boundary value PDEs \footnote{For more information about NDSolve, see Wolfram documentation center at https://reference.wolfram.com/language/ref/NDSolve.html}. We then calculate the expectation and variance of the solutions of our set of boundary value problems of PDEs. In Fig.\ref{fig:04}, the errors in the calculations of $E(u)$ and $V(u)$ using $m=3$ and Poly-Sinc (with orthonormal Legendre) and the references from the $100$ PDEs are presented. Using the spatial $L_2$-norm error, calculating the error in both $E(u)$ and $V(u)$ delivers error of order $\mathcal{O}(10^{-4})$ and $\mathcal{O}(10^{-6})$, respectively. In Fig.\ref{fig:05}, the error between the solution of the SPDE in (\ref{Eqn_100}), using the method in this paper, and the reference solution is presented. We choose four instances of $\xi$.  
\end{expt}

\begin{figure}[H]
\centering
\begin{subfigure}{.5\textwidth}
  \centering
  \includegraphics[width= .9\linewidth]{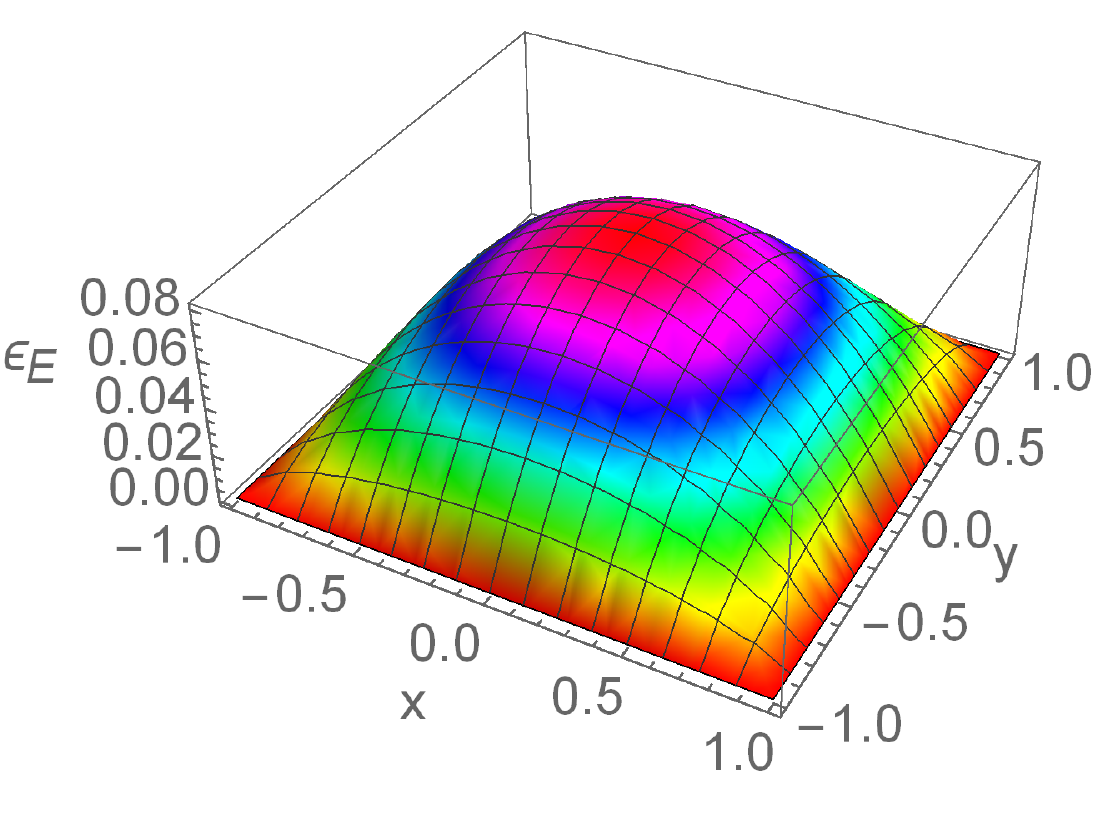}
  \caption{Absolute error in $E(u)$.}
\end{subfigure}%
\begin{subfigure}{.5\textwidth}
  \centering
  \includegraphics[width= .9\linewidth]{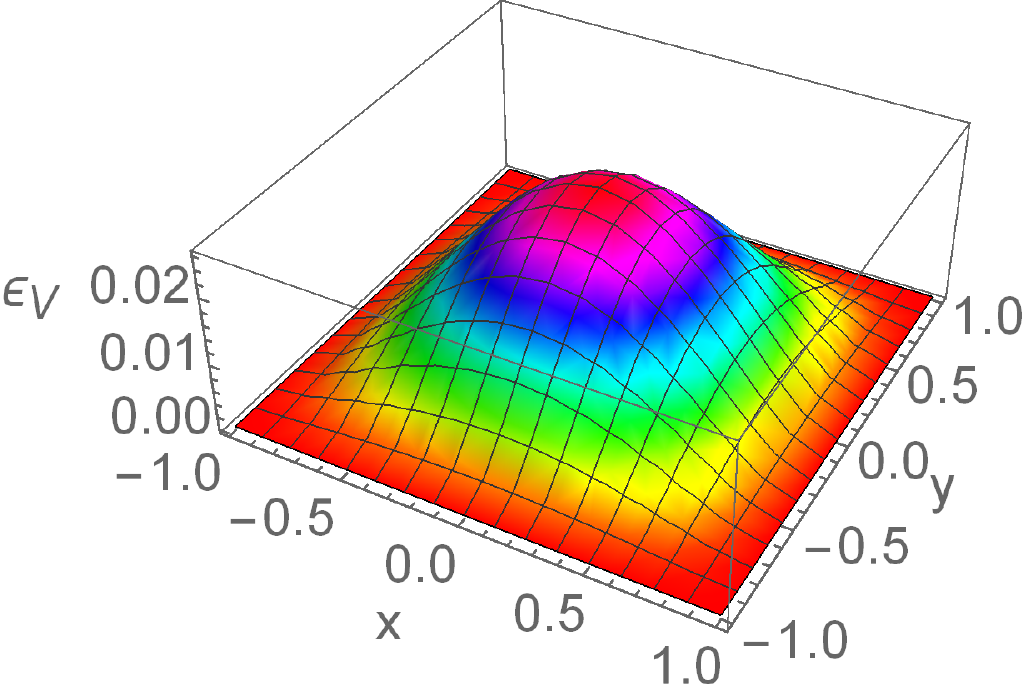}
  \caption{Absolute error in $V(u)$}
  \label{fig:sub2}
\end{subfigure}
\caption{Absolute error between the FD calculation and the calculations obtained from $100$ solutions.}
\label{fig:fd1}
\end{figure}

\begin{figure}[H]
\centering
\begin{subfigure}{.5\textwidth}
  \centering
  \includegraphics[width= .9\linewidth]{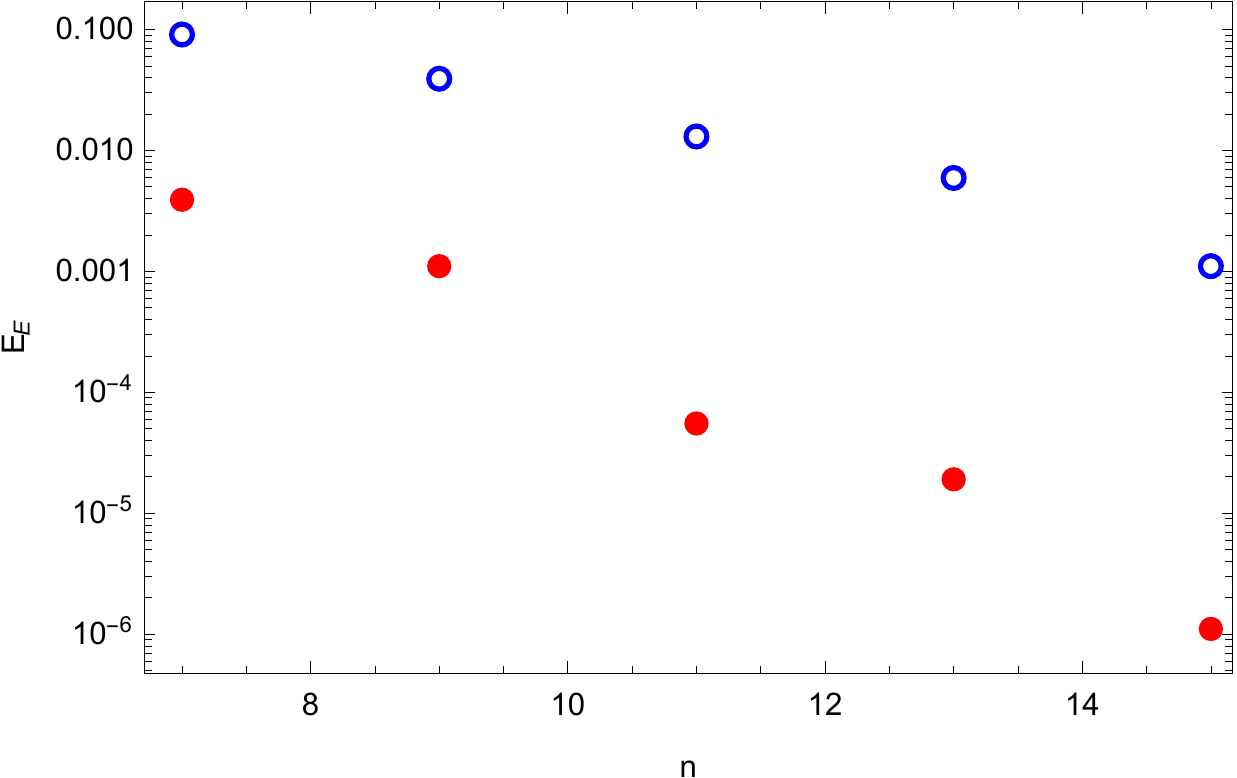}
  \caption{$L_2$ error in $E(u)$.}
\end{subfigure}%
\begin{subfigure}{.5\textwidth}
  \centering
  \includegraphics[width= .9\linewidth]{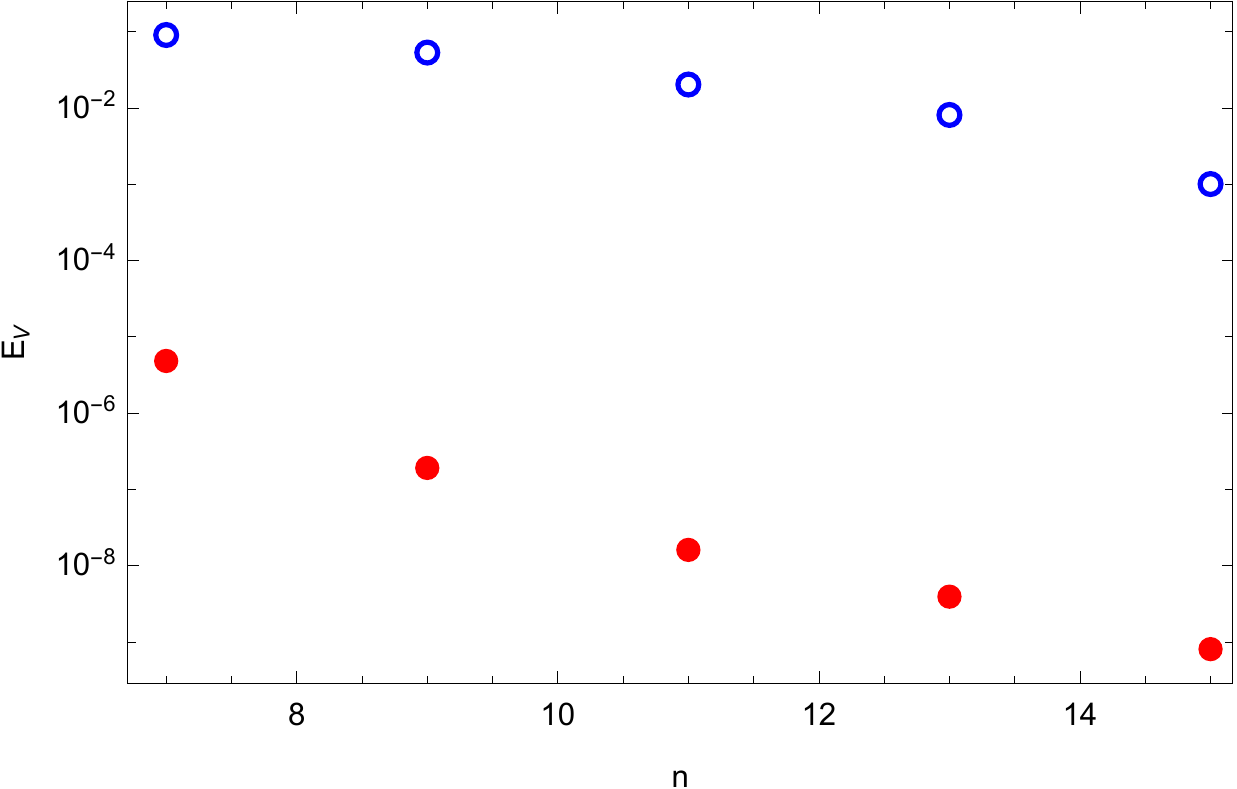}
  \caption{$L_2$ error in $V(u)$}
  \label{fig:sub2}
\end{subfigure}
\caption{Spatial $L_2$-error. The red dots for Poly-Sinc calculations and the blue circles for FD method with uniform meshes.}
\label{fig:fd2}
\end{figure}

\newpage
\begin{expt}{\bf{Comparison}}
\newline
In this experiment we compare the Poly-Sinc solution with the classical finite difference (FD) solution. In 5-point-star FD method \cite{Grossmann_2007}, we use an $11 \times 11$ meshing with constant step size for the spatial variables $x$ and $y$, which is the same number of Sinc points used in the Poly-Sinc solution. The error between finite difference solution and the reference exact solution is given in Fig.\ref{fig:fd1}. Using the spatial $L_2$-norm error, calculating the error in both $E(u)$ and $V(u)$ delivering error of order $\mathcal{O}(10^{-2})$. These calculations shows that for the same number of points, Poly-Sinc delivers better approximation for the solution of the SPDE. In Fig.\ref{fig:fd2} we run the calculations for different numbers of Sinc points $n=2N+1$ and use the same number of points in the FD method. We then calculate the $L_2$-norm error. Fig.\ref{fig:fd2} shows that the decaying rate of the error, in both mean and variance, is better in Poly-Sinc than the FD method. Moreover, the Poly-Sinc decaying rates of errors are following qualitatively the upper bound in formula (\ref{eq:coler}).
\end{expt}

\newpage
\subsection{Multiple Stochastic Variables}
We solve the model problem defined in Section \ref{sec:2} for five stochastic variables, cf. \cite{GITTELSON_2013}. Consider the SPDE defined in (\ref{Eqn_01}) with $K=5$ in (\ref{Eqn_02}) and where,

\begin{align*}
a_1(x,y) & =  \mathsmaller{\frac{1}{4}} \cos(2 \pi x)\\
a_2(x,y) & =  \mathsmaller{\frac{1}{4}} \cos(2 \pi y)\\
a_3(x,y) & =  \mathsmaller{\frac{1}{16}} \cos(4 \pi x)\\
a_4(x,y) & =  \mathsmaller{\frac{1}{16}} \cos(4 \pi y)\\
a_5(x,y) & =  \mathsmaller{\frac{1}{8}} \cos(2 \pi x) \cos(2 \pi y).
\end{align*}

\(\Theta =\left\{\xi_k\right\}^{5}_{k=1}\) is a set of independent random variables uniformly distributed in $[-1,1]$. For this SPDE we run four experiments.

\begin{figure}[H]
\centering
\begin{subfigure}{.5\textwidth}
  \centering
  \includegraphics[width= \linewidth]{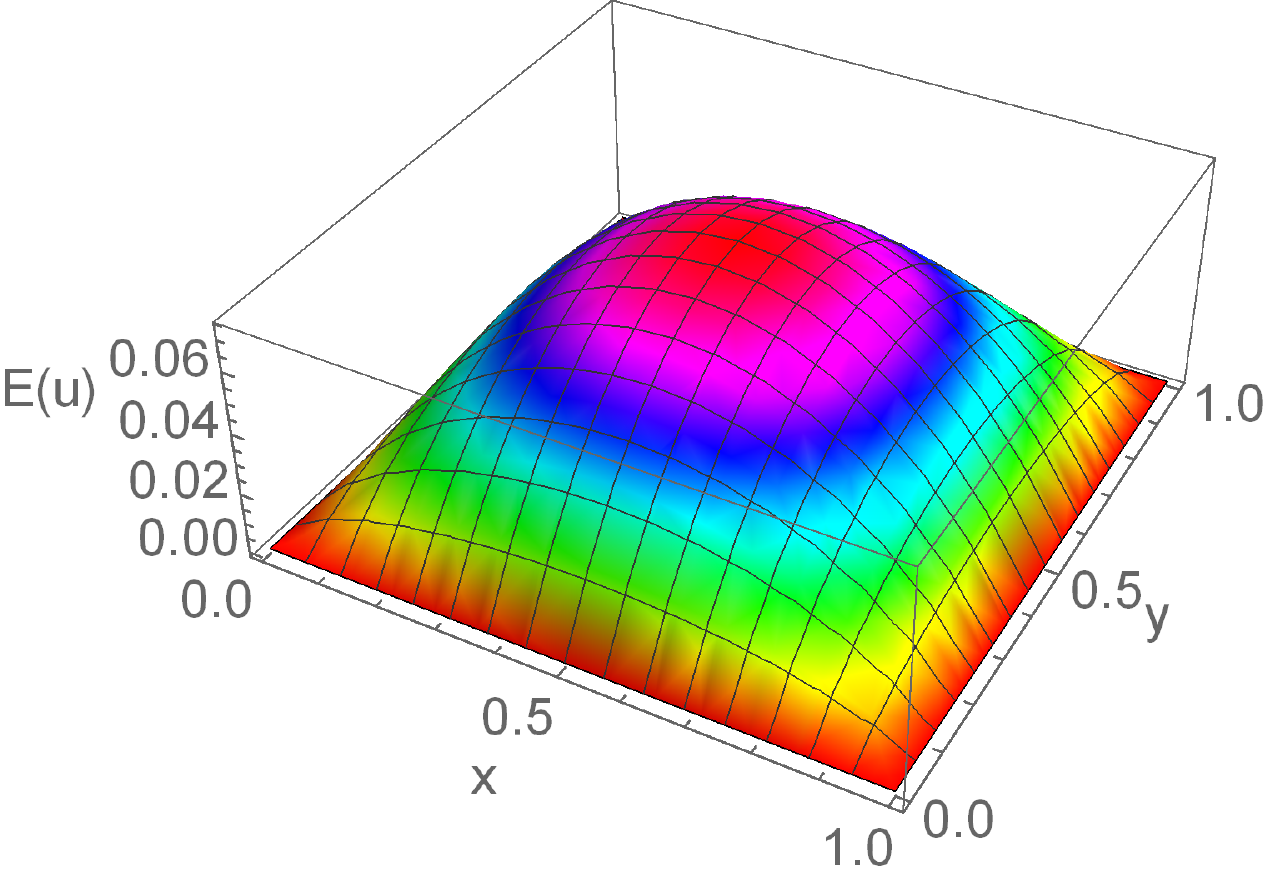}
  \caption{$E(u)$.}
\end{subfigure}%
\begin{subfigure}{.5\textwidth}
  \centering
  \includegraphics[width=.9\linewidth]{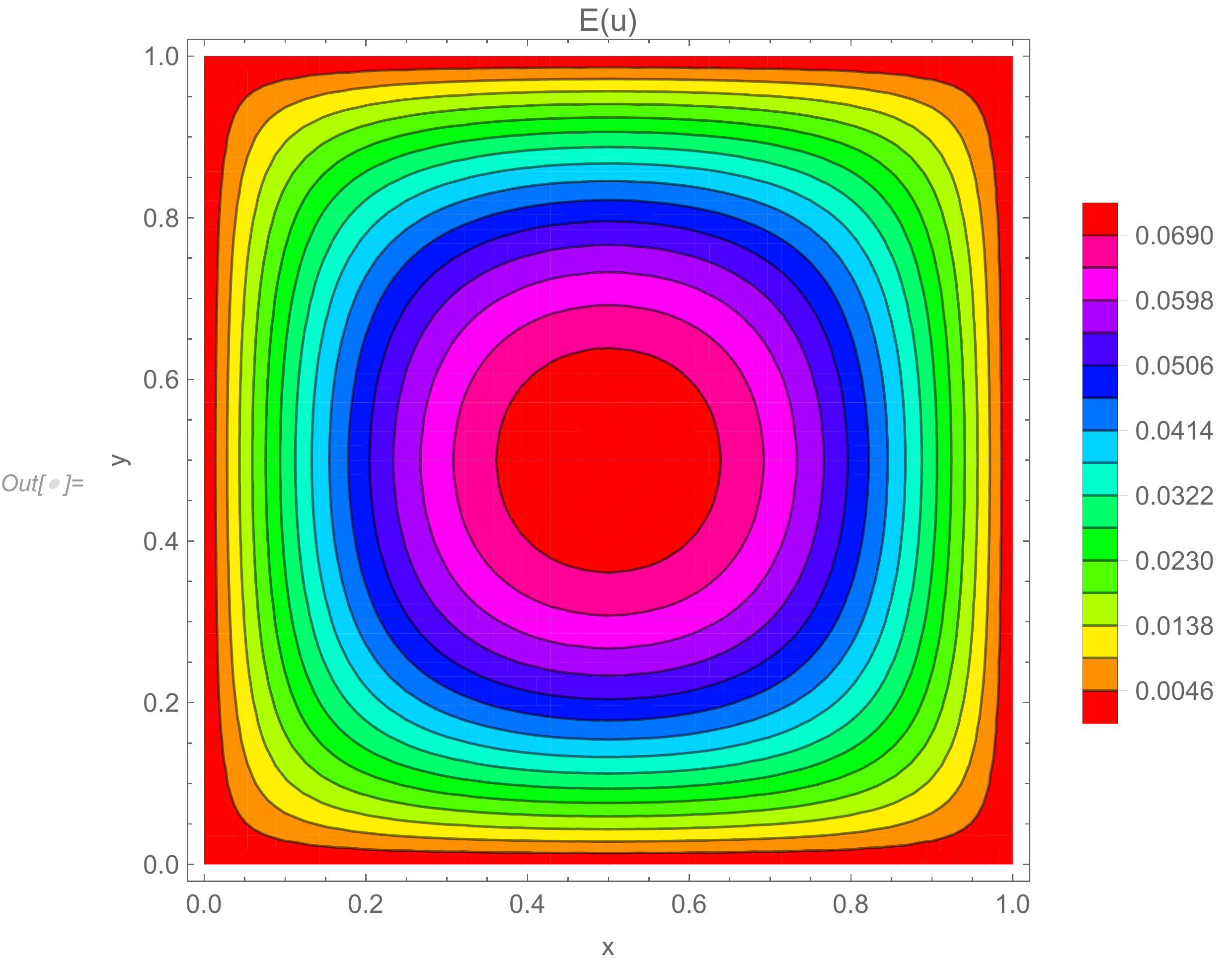}
  \caption{Contour plot of $E(u)$}
  \label{fig:sub2}
\end{subfigure}
\caption{The expectation, $E(u)$, using $K=5,\, P=3$ and Poly-Sinc with $N=5$.}
\label{fig:08}
\end{figure}

\begin{figure}[H]
\centering
\begin{subfigure}{.5\textwidth}
  \centering
  \includegraphics[width= \linewidth]{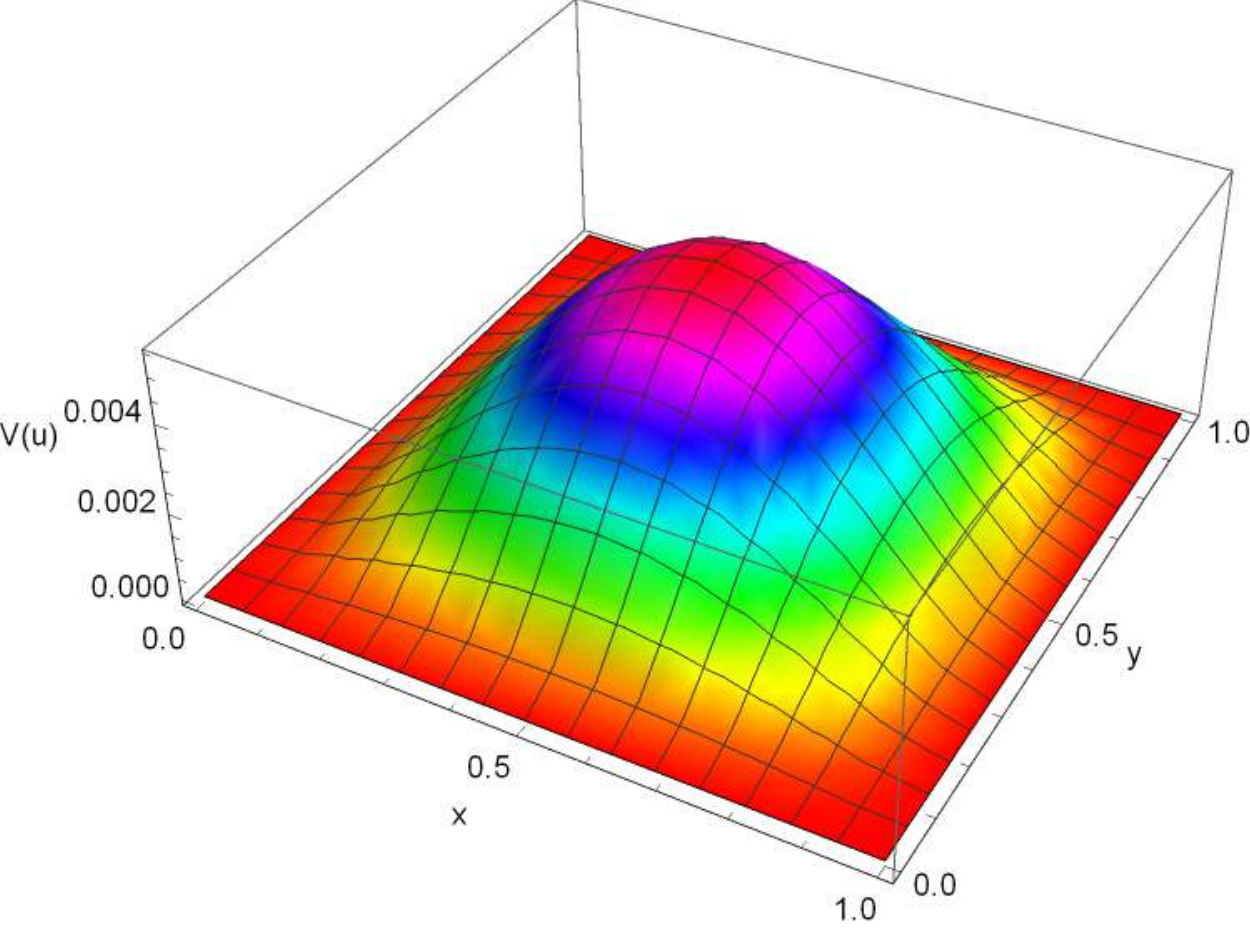}
  \caption{$V(u)$.}
\end{subfigure}%
\begin{subfigure}{.5\textwidth}
  \centering
  \includegraphics[width=.9\linewidth]{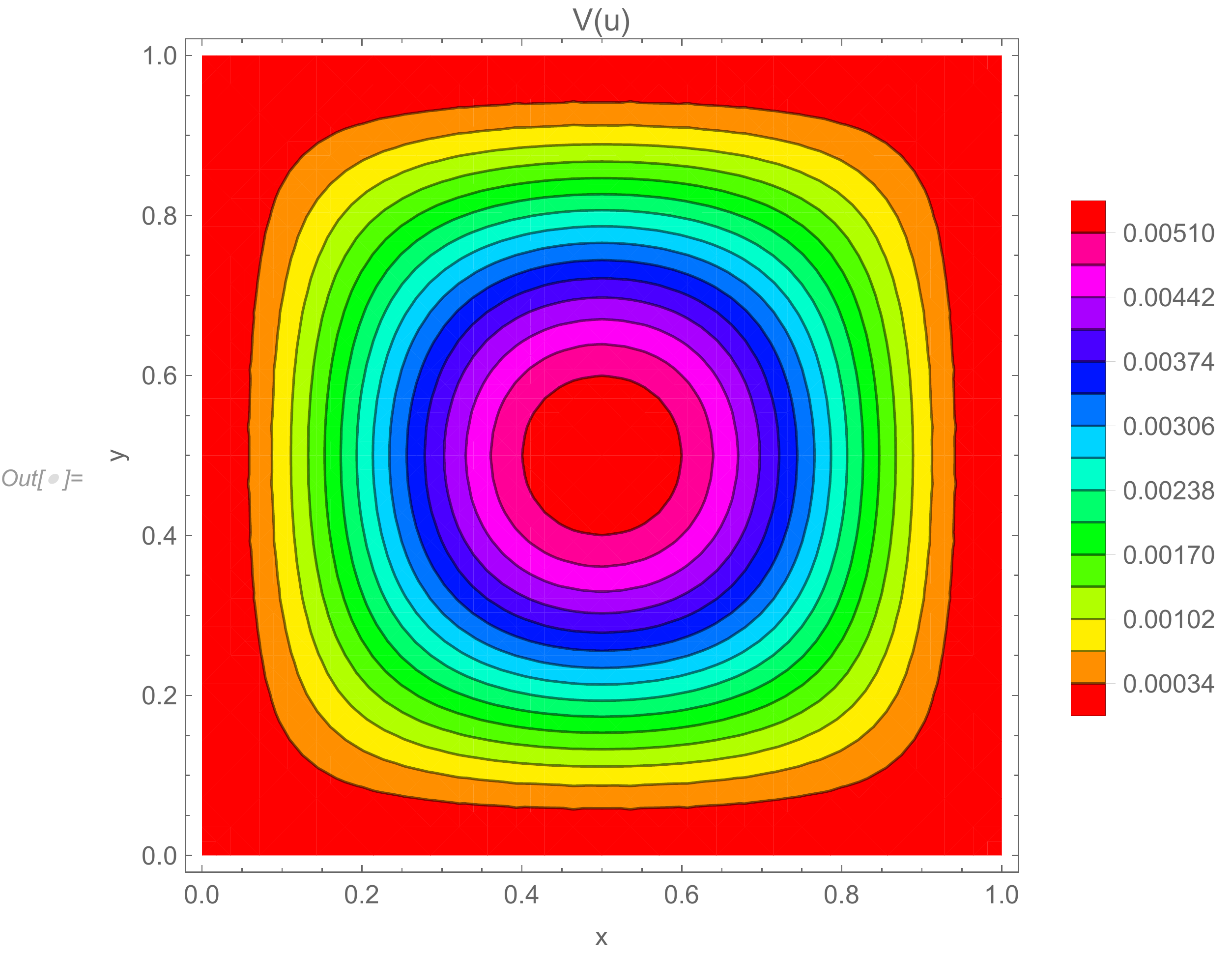}
  \caption{Contour plot of $V(u)$}
  \label{fig:sub2}
\end{subfigure}
\caption{The variance, $V(u)$, using $K=5,\, P=3$ and Poly-Sinc with $N=5$.}
\label{fig:09}
\end{figure}

\begin{expt}{$\mathbf{E(u)}$ \bf{and} $\mathbf{V(u)}$}
\newline
In this experiment, we perform the Galerkin method along side the multivariate PC. For the PC parameters, we choose $K=5$ and $P=3$. Due to (\ref{Eqn_09}), the number of multivariate Legendre polynomials is $m+1=56$. As a result the three-dimensional array $\left\langle  \xi _k \Phi _i(\Theta ),\Phi _j(\Theta )\right\rangle$ is of dimension $5 \times 56\times 56$. For the Poly-Sinc solution of the resulting system of PDEs, we use $N=5$, i.e. $n=11$ Sinc points. In Fig. \ref{fig:08} and Fig.\ref{fig:09} the expectation and variance plots are presented.
\end{expt}

\begin{figure}[H]
\centering
\includegraphics[scale=0.7]{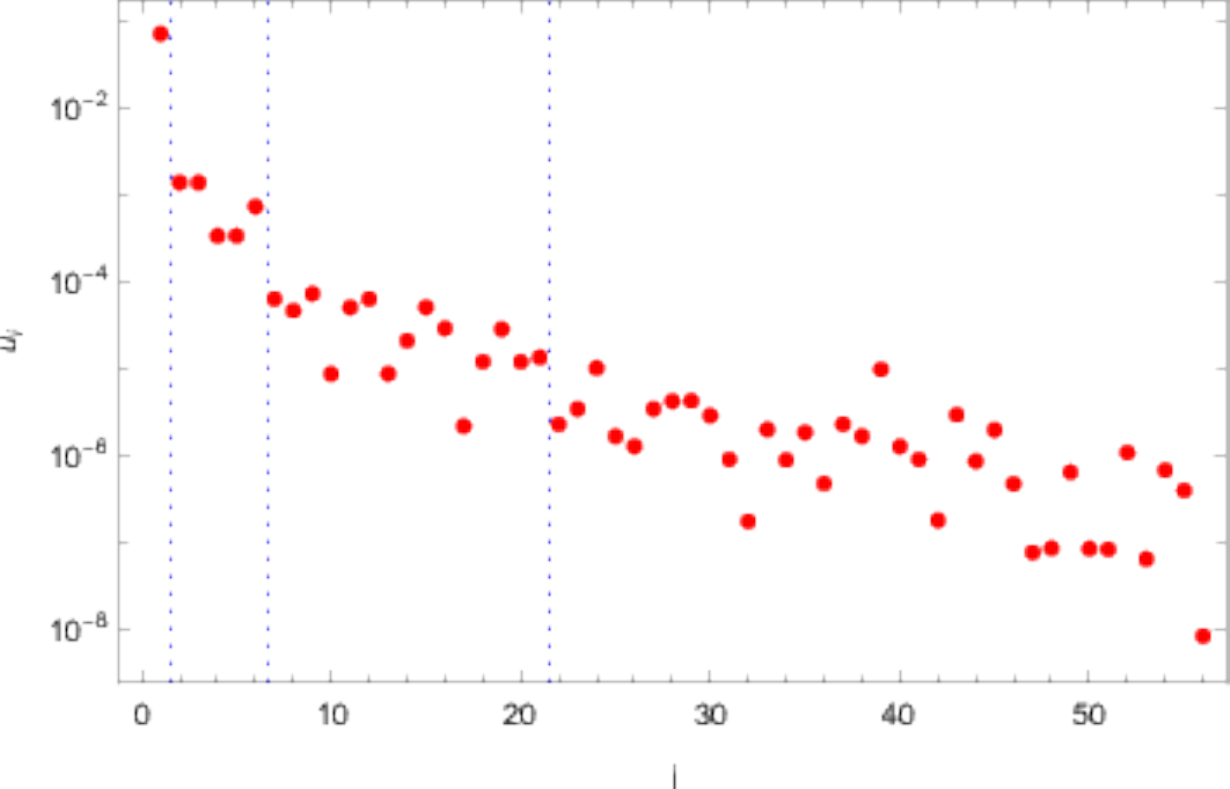}
\caption{Logarithmic plot of maximum of coefficient functions $u_{i-1}(x,y)$ for $i=1,\ldots,56$. The dotted lines separate the degrees of basis polynomials.}
\label{fig:10}
\end{figure}

\begin{expt}{\bf{Coefficients Functions}}
\newline
Similar to the second experiment in Example 1, we would like to study the accuracy of the polynomial expansion. In other words, study the decaying rate, to zero, of these functions. In Fig.\ref{fig:10}, the first six coefficients functions of the Poly-Sinc solution are given. These six coefficient functions are associated to the basis polynomials of degree zero and one. In Fig. \ref{fig:10}, the logarithmic plot of the maximum of the absolute value of the coefficient functions $u_{i-1}(x,y),\,\,i=1,\ldots,56$ on the spatial domain is presented. We can see the fast decaying rate to zero.
\end{expt}

\begin{figure}[H]
\centering
\includegraphics[width= \linewidth]{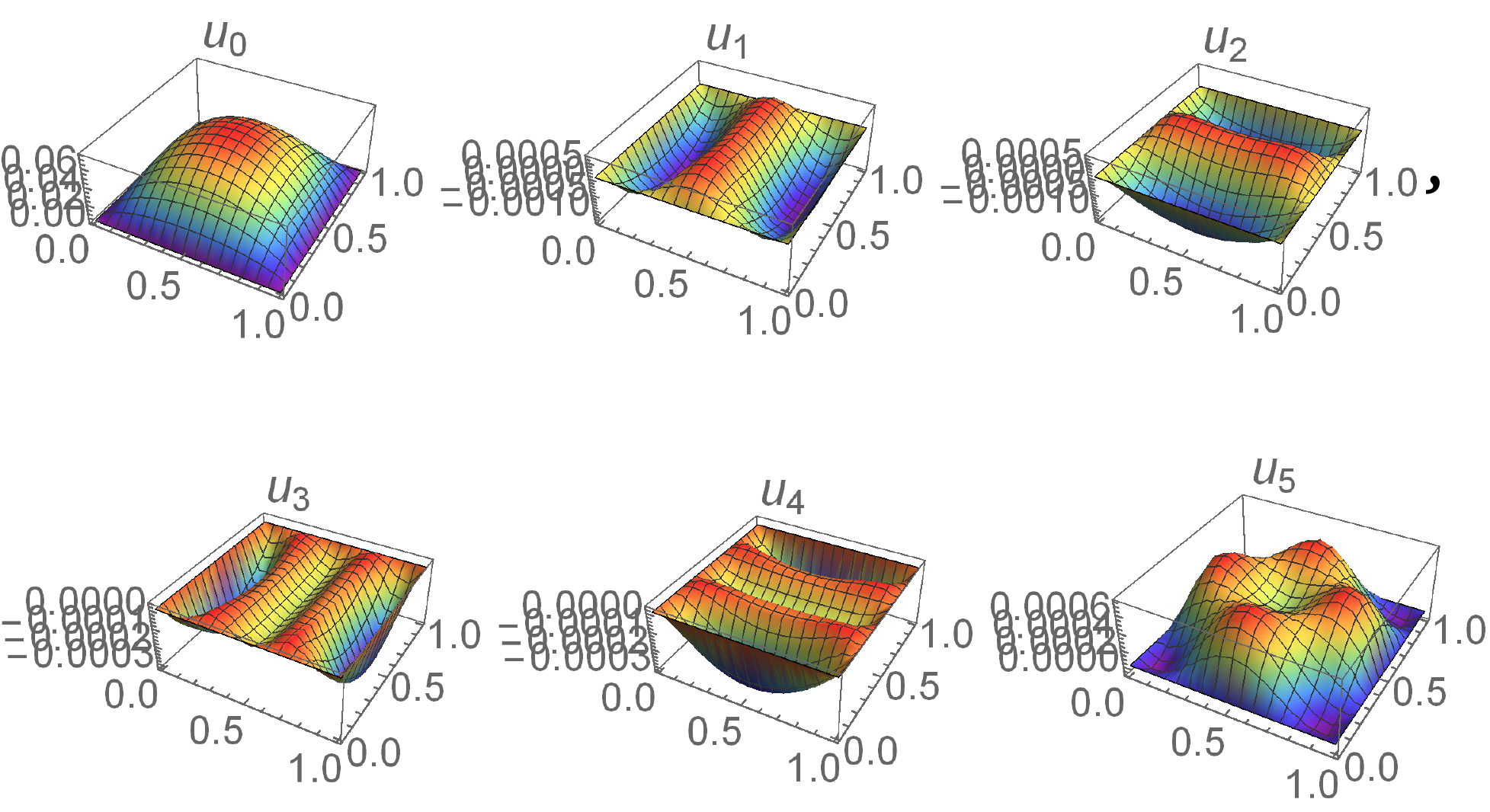}
\caption{Coefficients functions $u_{i}(x,y),\,\,i=0,1,\ldots,5$.}
\label{fig:1001}
\end{figure}

\begin{figure}[H]
\centering
\begin{subfigure}{.5\textwidth}
  \centering
  \includegraphics[width= .9\linewidth]{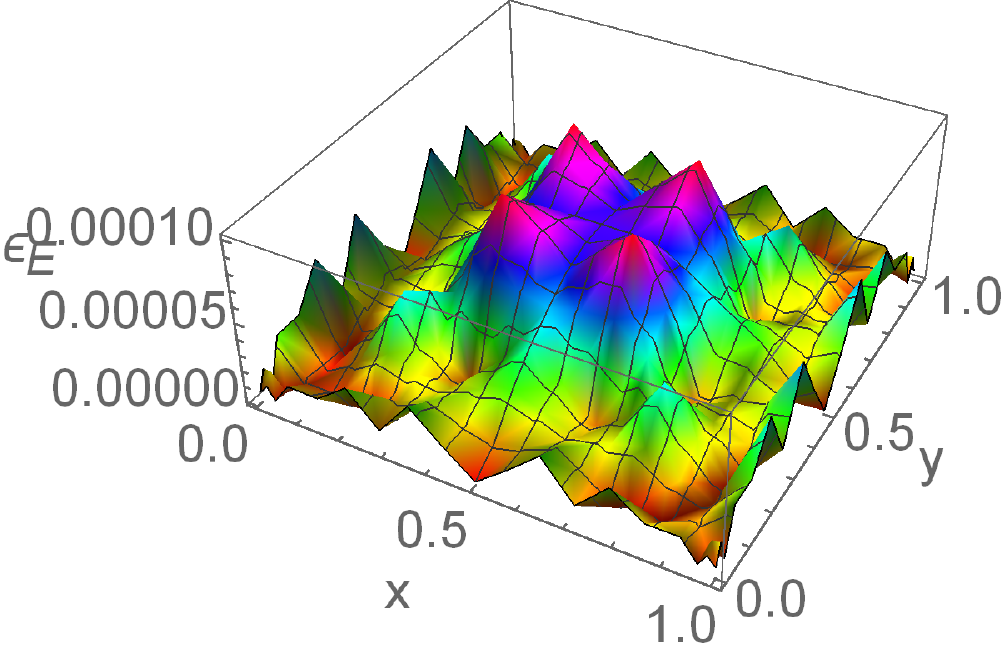}
  \caption{Absolute error in $E(u)$.}
\end{subfigure}%
\begin{subfigure}{.5\textwidth}
  \centering
  \includegraphics[width=.9\linewidth]{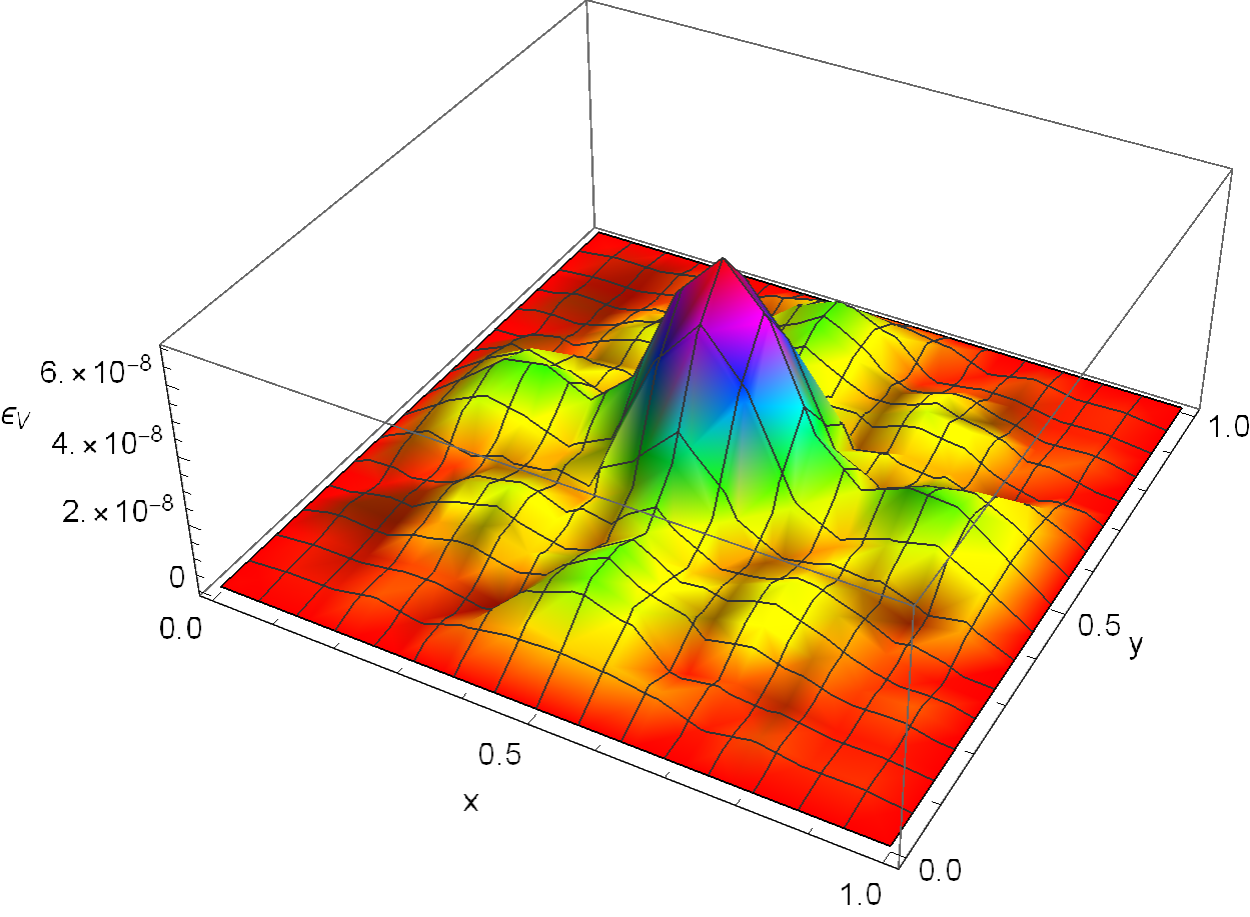}
  \caption{Absolute error in $V(u)$}
  \label{fig:sub2}
\end{subfigure}
\caption{Absolute error between the Poly-Sinc calculation and the FE.}
\label{fig:11}
\end{figure}

\begin{expt}{\bf{Error}}
\newline
The idea of creating a set of (exact) instance solutions we used in the previous example is not applicable here as we have a set of $5$ random variables. For that we need to find a different reference to check the accuracy of our solution. We use the Finite Element (FE) solution with cell meshing $10^{-3}$ to solve the resulting system of PDEs. The FE element method is a part of the package NDSolve"FEM" in Mathematica 11 that uses the rectangular meshing of the domain and Dirichlet boundary conditions \footnote{For more information about NDSolve "FEM", see Wolfram documentation center at https://reference.wolfram.com/language/FEMDocumentation/guide/FiniteElementMethodGuide.html}. In Fig.\ref{fig:11}, the error for the expectation and variance is presented. Using the $L_2$-norm error, calculating the error in both $E(u)$ and $V(u)$ deliver error of order $\mathcal{O}(10^{-4})$ and $\mathcal{O}(10^{-8})$, respectively. 
\end{expt}

\begin{figure}[H]
\centering
\begin{subfigure}{.5\textwidth}
  \centering
  \includegraphics[width= 0.95\linewidth]{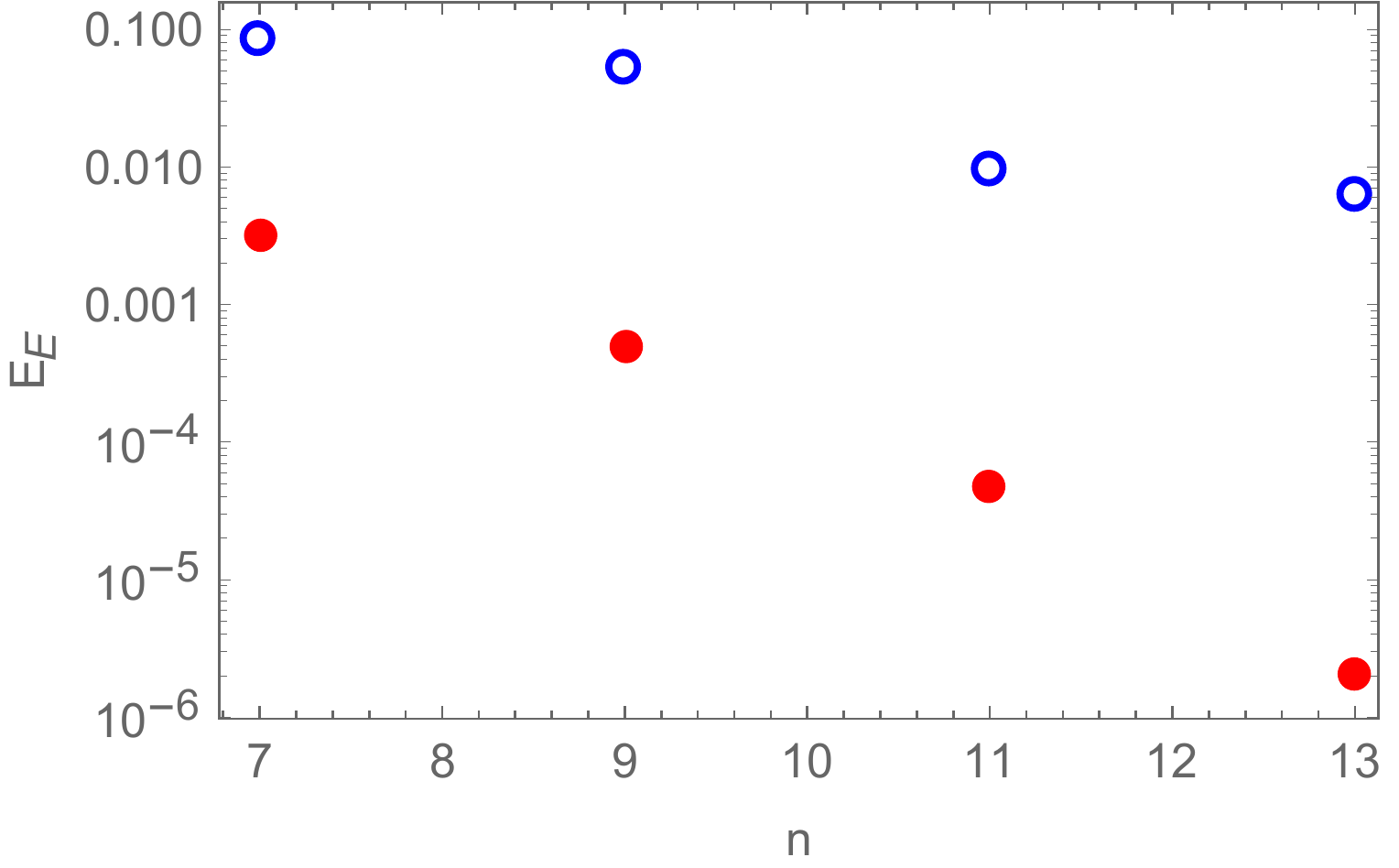}
  \caption{$L_2$ error in $E(u)$.}
\end{subfigure}%
\begin{subfigure}{.5\textwidth}
  \centering
  \includegraphics[width=.9\linewidth]{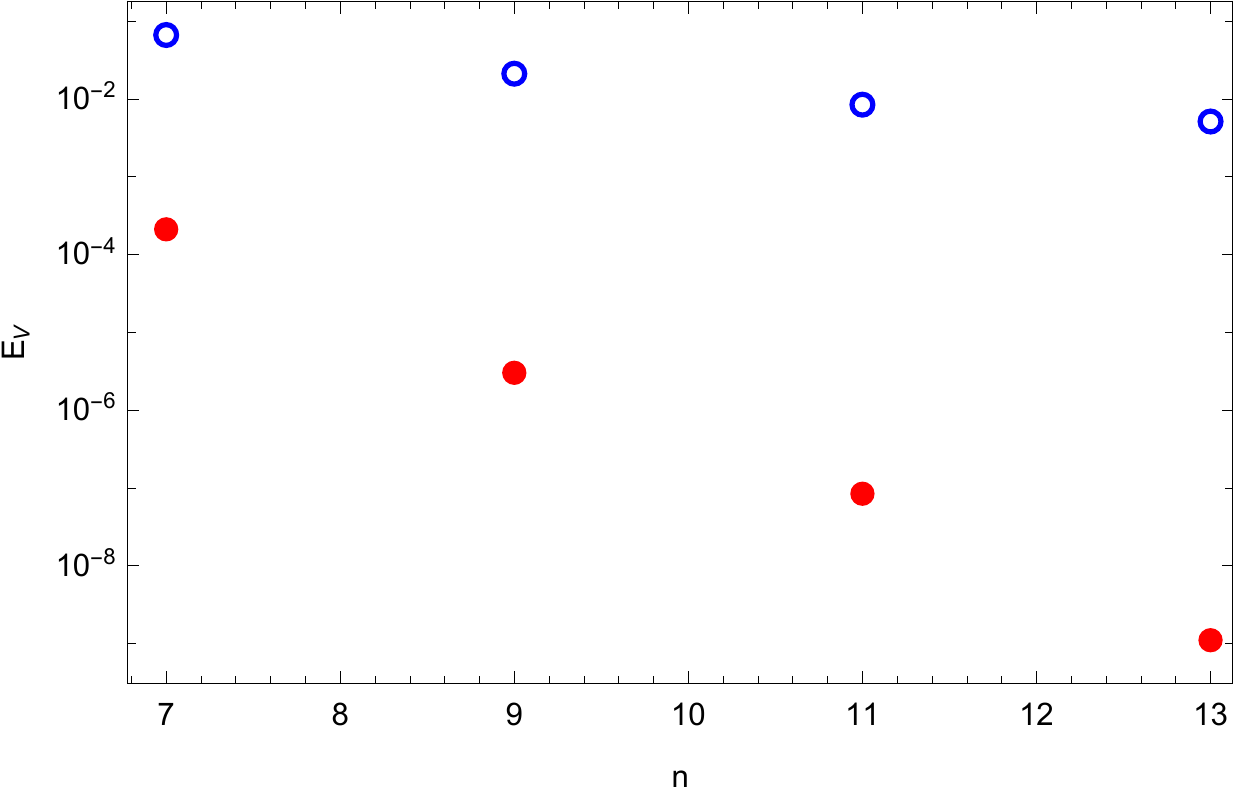}
  \caption{$L_2$ error in $V(u)$}
  \label{fig:sub2}
\end{subfigure}
\caption{Spatial $L_2$-error. The red dots for Poly-Sinc calculations and the blue circles for FD method.}
\label{fig:fd12}
\end{figure}

\newpage
\begin{expt}{\bf{Comparison}}
\newline
In this experiment we compare the Poly-Sinc solution with the 5-point-star FD method. The reference solution is the Finite Element (FE) solution with cell meshing $10^{-3}$. In Fig.\ref{fig:fd12} we run the calculations for different numbers of Sinc points $n=2N+1$ and use the same number of points in FD. We then calculate the $L_2$-norm error. These calculations show that the decaying rate of the error, in both mean and variance, is better in Poly-Sinc than the FD method. Moreover, the Poly-Sinc decaying rates of errors are following qualitatively the exponential decaying rate in (\ref{eq:coler}).

\end{expt}

\section{Conclusion}
In this work we have formulated an efficient and accurate collocation scheme for solving a system of elliptic PDEs resulting from an SPDE. The idea of the scheme is to use a small number of collocation points to solve a large system of PDEs. We introduced the collocation theorem based on the error rate and the Lebesgue constant of the 2D Poly-Sinc approximation. As applications, we discussed two examples, the first example with one random variable while the other with five random variables. For each case the expectation, variance, and error are discussed. The experiments show that using Poly-Sinc approximation to solve the system of PDEs is an efficient method. The number of Sinc points needed to get this accuracy is small and the error decays faster than in the classical techniques, as the finite difference method.



\end{document}